\documentclass[reqno,11pt]{article}
\usepackage{amsmath, amstexnb, amsthm}
\usepackage{amssymb}  

\input{geom.sty}

\begin{document}


\title{Geometric singular perturbation theory \\
for stochastic differential equations}
\author{Nils Berglund and Barbara Gentz}
\date{}   

\maketitle

\begin{abstract}
\noindent
We consider slow--fast systems of differential equations, in which both the
slow and fast variables are perturbed by noise. When the deterministic system
admits a uniformly asymptotically stable slow manifold, we show that the
sample paths of the stochastic system are concentrated in a neighbourhood of
the slow manifold, which we construct explicitly. Depending on the dynamics of
the reduced system, the results cover time spans which can be exponentially
long in the noise intensity squared (that is, up to Kramers' time). We obtain
exponentially small upper and lower bounds on the probability of exceptional
paths.  
If the slow manifold contains bifurcation points, we show similar
concentration properties for the fast variables corresponding to
non-bifurcating modes. We also give conditions under which the system can
be approximated by a lower-dimensional one, in which the fast variables
contain only bifurcating modes.  
\end{abstract}

\leftline{\small{\it Date.\/} March 30, 2002. {\it Revised.\/} October 28,
  2002.}
\leftline{\small 2000 {\it Mathematical Subject Classification.\/} 
37H20, 34E15 (primary), 60H10 (secondary)
}

\noindent{\small{\it Running head.\/}
Geometric singular perturbation theory for SDEs.}


\section{Introduction}

Systems involving two well-separated timescales are often described by
slow--fast differential equations of the form
\begin{equation}   \label{i1}
\begin{split}
\eps\dot x &= f(x,y,\eps), \\
\dot y &= g(x,y,\eps),
\end{split}
\end{equation}
where $\eps$ is a small parameter. Since $\dot x$ can be much larger than
$\dot y$, $x$ is called the \defwd{fast variable} and $y$ is called the
\defwd{slow variable}. Such equations occur, for instance, in climatology,
with the slow variables describing the state of the oceans, and the fast
variables the state of the atmosphere. In physics, slow--fast equations 
model in particular systems containing heavy particles (e.\,g.\ nuclei) and
light particles (e.\,g.\ electrons). Another example, taken from ecology,
would be the dynamics of a predator--prey system in which the rates of
reproduction of predator and prey are very different.  

The system~\eqref{i1} behaves singularly in the limit $\eps\to0$. In fact,
the results depend on the way this limit is performed. If we simply set
$\eps$ to zero in~\eqref{i1}, we obtain the algebraic--differential system 
\begin{equation}   \label{i2}
\begin{split}
0 &= f(x,y,0), \\
\dot y &= g(x,y,0).
\end{split}
\end{equation}
Assume there exists a differentiable manifold with equation $x=x^\star(y)$
on which $f=0$. Then $x=x^\star(y)$ is called a \defwd{slow manifold}, and
the dynamics on it is described by the \defwd{reduced equation} 
\begin{equation}   \label{i3}
\dot y = g(x^\star(y),y,0). 
\end{equation}
Another way to analyze the limit $\eps\to0$ is to scale time by a factor
$1/\eps$, so that the slow--fast system~\eqref{i1} becomes
\begin{equation}   \label{i4}
\begin{split}
x^\prime &= f(x,y,\eps), \\
y^\prime &= \eps g(x,y,\eps). 
\end{split}
\end{equation}
In the limit $\eps\to0$, we obtain the so-called \defwd{associated system} 
\begin{equation}   \label{i5}
\begin{split}
x^\prime &= f(x,y,0), \\
y^\prime &= 0,
\end{split}
\end{equation}
in which $y$ plays the r\^ole of a parameter. The slow manifold
$x=x^\star(y)$ consists of equilibrium points of~\eqref{i5}, and~\eqref{i4}
can be viewed as a perturbation of~\eqref{i5} with slowly drifting parameter
$y$. 

Under certain conditions, both the reduced equation~\eqref{i3} and the
associated system~\eqref{i5} give good approximations of the initial
slow--fast system~\eqref{i1}, but on different timescales. Assume for
instance that for each $y$, $x^\star(y)$ is an asymptotically stable
equilibrium of the associated system~\eqref{i5}. Then solutions
of~\eqref{i1} starting in a neighbourhood of the slow manifold will approach
$x^\star(y)$ in a time of order $\eps\abs{\log\eps}$. During this time
interval  they are well approximated by solutions of~\eqref{i5}. This first
phase of the motion is sometimes called the \defwd{boundary-layer}
behaviour. For larger times, solutions of~\eqref{i1} remain in an
$\eps$-neighbourhood of the slow manifold, and are thus well approximated by
solutions of the reduced equation~\eqref{i3}. This result was first proved
by Grad{\v{s}}te{\u\i}n \cite{Grad} and Tihonov \cite{Tihonov}. 

Fenichel \cite{Fenichel} has given results allowing for a geometrical
description of these phenomena in terms of invariant manifolds. He showed,
in particular, the existence of an invariant manifold
\begin{equation}   \label{i6}
x = \bx(y,\eps),
\qquad\qquad
\text{with $\bx(y,\eps) = x^\star(y) + \Order{\eps}$,}
\end{equation}
for sufficiently small $\eps$, whenever $x^\star(y)$ is a family of
hyperbolic equilibria of the associated system~\eqref{i5}. The dynamics on
this invariant manifold is given by the equation 
\begin{equation}   \label{i7}
\dot y = g(\bx(y,\eps),y,\eps),
\end{equation}
which can be treated by methods of regular perturbation theory, and reduces
to~\eqref{i3} in the limit $\eps\to0$. In fact, Fenichel's results are more
general. For instance, if $x^\star(y)$ is a saddle, they also show the
existence of invariant manifolds associated with the stable and unstable
manifolds of $x^\star(y)$. See \cite{Jones} for a review. 

New, interesting phenomena arise when the dynamics of~\eqref{i7} causes $y$
to approach a bifurcation point of~\eqref{i5}. For instance, the passage
through a saddle--node bifurcation, corresponding to a fold of the slow
manifold, produces a jump to some other region in phase space, which can
cause relaxation oscillations and hysteresis phenomena (see in particular
\cite{Pontryagin} and \cite{Haberman}, as well as \cite{MishchenkoRozov} for an
overview). Transcritical and pitchfork bifurcations generically lead to a
smoother transition to another equilibrium
\cite{LebovitzSchaar1,LebovitzSchaar2}, while the passage through a Hopf
bifurcation is accompanied by the delayed appearance of oscillations
\cite{Neishtadt1,Neishtadt2}. There exist many more recent studies of what
has become known as the field of dynamic bifurcations, see for instance
\cite{Benoit}. 

In many situations, low-dimensional ordinary differential equations of the
form $\dot x= f(x)$ are not sufficient to describe the dynamics of the
system under study. The effect of unknown degrees of freedom is often
modelled by noise, leading to a stochastic differential equation (SDE) of
the form 
\begin{equation}   \label{i8}
\6x_t = f(x_t)\6t + \sigma F(x_t) \6W_t,
\end{equation}
where $\sigma$ is a small parameter, and $W_t$ denotes a standard, generally
vector-valued Brownian motion. On short timescales, the main effect of
the noise term $\sigma F(x_t) \6W_t$ is to cause solutions to fluctuate
around their deterministic counterpart, but the probability of large
deviations is very small (of the order $\e^{-{\it const}/\sigma^2}$). On
longer timescales, however, the noise term can induce transitions to other
regions of phase space. 

The best understood situation is the one where $f$ admits an asymptotically
stable equilibrium point $x^\star$. The first-exit time $\tau(\w)$ of the
sample path $x_t(\w)$ from a neighbourhood of $x^\star$ is a random variable,
the characterization of which is the object of the \defwd{exit problem}. If $f$
derives from a potential $U$ (i.\,e., $f=-\nabla U$) of which $x^\star$ is a
local minimum, the asymptotic behaviour of the typical first-exit time for
$\sigma\ll 1$ has been long known by physicists: it is of order
$\e^{2H/\sigma^2}$, where $H$ is the height of the lowest potential
barrier separating $x^\star$ from other potential wells. A theory of large
deviations generalizing this result to quite a large class of SDEs has been
developed by Freidlin and Wentzell~\cite{FW}. More detailed information on
the asymptotics of the expected first-exit time has been obtained,
see~\cite{Azencott,FJ} and the very precise results by Bovier, Eckhoff,
Gayrard and Klein~\cite{BEGK,BGK} on the relation between the expected
first-exit time, capacities and the spectrum of the generator of the
diffusion. The distribution of $\tau$ has been studied by Day~\cite{Day1}.

The more difficult problem of the dynamics near a saddle point has been
considered in \cite{Kifer} and in \cite{Day2}. The situation where $f$ depends
on a parameter and undergoes bifurcations has not yet been studied in that
much detail. An approach based on the notion of random attractors
\cite{Schmal,Arnold1,CF1} gives information on the limit $t\to\infty$, when
the system has reached a stationary state. Note, however, that the time
needed to reach this regime, in which (in the gradient case) $x_t$ is most
likely to be found near the deepest potential well, may be very long if the
wells are separated by barriers substantially higher than $\sigma^2$. The
dynamics on intermediate timescales, known as the \defwd{metastable regime},
is not yet well understood in the presence of bifurcations. 

In this work, we are interested in the effect of noise on
slow--fast systems of the form~\eqref{i1}. Such systems have been studied
before in~\cite{Freidlin2}, using techniques from large deviation theory to
describe the limit $\sigma\to0$. Here we use different methods to give a 
more precise description of the regime of small, but finite noise
intensity, our main goal being to estimate quantitatively the noise-induced
spreading of typical paths, as well as the probability of exceptional
paths.  We will consider situations in which both the slow and fast
variables are affected by noise, with noise intensities taking into account
the difference between the timescales. In~\eqref{i8}, the diffusive nature
of the Brownian motion causes paths to spread like $\sigma\sqrt t$. In the
case of the slow--fast system~\eqref{i1}, we shall choose the following
scaling of the noise intensities:
\begin{equation}   \label{i9}
\begin{split}
\6x_t &= \frac1\eps f(x_t,y_t,\eps)\6t + \frac\sigma{\sqrt\eps}
F(x_t,y_t,\eps)\6W_t, \\
\6y_t &= g(x_t,y_t,\eps)\6t + \sigma^\prime G(x_t,y_t,\eps) \6W_t.
\end{split}
\end{equation}
In this way, $\sigma^2$ and $(\sigma^\prime)^2$ both measure the ratio
between the rate of diffusion squared and the speed of drift,
respectively,  for the fast and slow variable. We consider general
finite-dimensional $x\in\R^n$ and $y\in\R^m$, while $W_t$ denotes a
$k$-dimensional standard Brownian motion. Accordingly, $F$ and $G$ are
matrix-valued functions of respective dimensions $n~\times~k$ and $m\times
k$. The matrices $F(x,y,\eps)$ will be assumed to satisfy some (rather weak)
nondegeneracy condition, see Remark~\ref{rem_randst1}.
We consider $\eps$, $\sigma$ and $\sigma^\prime$ as small parameters,
and think of $\sigma$ and $\sigma^\prime$ as functions of $\eps$. We limit
the analysis to situations where $\sigma^\prime$ does not dominate $\sigma$,
i.\,e., we assume $\sigma^\prime=\rho\sigma$ where $\rho$ may depend on
$\eps$ but is uniformly bounded above in $\eps$.  

We first consider the case where the deterministic slow--fast
system~\eqref{i1} admits an asymptotically stable slow manifold
$x^\star(y)$. Our first main result, Theorem~\ref{thm1}, states that the
sample paths of ~\eqref{i9} are concentrated in a \lq\lq layer\rq\rq\
surrounding the adiabatic manifold $\bx(y,\eps)$, of the form 
\begin{equation}   \label{i10}
\cB(h) = 
\bigsetsuch{(x,y)}{\bigpscal{(x-\bx(y,\eps))}{\Xbar(y,\eps)^{-1} 
    (x-\bx(y,\eps))} < h^2}
\end{equation}
up to time $t$, with a probability behaving roughly like
$1-(t^2/\eps)\e^{-h^2/2\sigma^2}$ as long as the paths do not reach the
vicinity of a bifurcation point. The matrix $\Xbar(y,\eps)$, defining the
ellipsoidal cross-section of the layer, is itself a solution of a slow--fast
system, and depends only on the values of $F$ and $\sdpar fx$ on the slow
manifold. In particular, $\Xbar(y,0)$ is a solution of the Lyapunov
equation 
\begin{equation}   \label{i11}
A^\star(y)X + X\transpose{A^\star(y)} + 
F(x^\star(y),y,0) \transpose{F(x^\star(y),y,0)} = 0,
\end{equation}
where $A^\star(y) = \sdpar fx(x^\star(y),y,0)$. For instance, if $f$
derives from a potential $U$, $-A^\star$ is the Hessian matrix of $U$ at
its minimum, and $\cB(h)$ is more elongated in those directions in which
the curvature of $U$ is smallest. 

Theorem~\ref{thm2} gives a more detailed description of the dynamics inside
$\cB(h)$, by showing that paths $(x_t,y_t)$ are concentrated in a
neighbourhood of the deterministic solution $(\xdet_t,\ydet_t)$ at least up to
times of order $1$. The spreading in the $y$-direction grows at a rate
corresponding to the finite-time Lyapunov exponents of the deterministic
solution.  

Next we turn to situations where the deterministic solution approaches a
bifurcation point of the associated system. In this case, the adiabatic
manifold $\bx(y,\eps)$ is not defined in general. However, by splitting $x$
into a stable direction $x^-$ and a bifurcating direction $z$, one can
define a (centre) manifold $x^- = \bxm(z,y,\eps)$ which is locally
invariant under the deterministic flow. Theorem~\ref{thm3} shows that paths
of the stochastic system are concentrated in a neighbourhood of
$\bxm(z,y,\eps)$. The  size of this neighbourhood again depends on noise and
linearized drift term in the stable $x^-$-direction. 

In order to make use of previous results on the passage through bifurcation
points for one-dimensional fast variables, such as \cite{BG1,BG2,BG3}, it
is necessary to control the deviation between solutions of the full
system~\eqref{i9}, and the reduced stochastic system obtained by setting
$x^-$ equal to $\bxm(z,y,\eps)$. Theorem~\ref{thm4} provides such an
estimate under certain assumptions on the dynamics of the reduced system. 

We present the detailed results in Section~\ref{sec_results},
Subsection~\ref{ssec_detst} containing a summary of results on
deterministic slow--fast systems, while Subsection~\ref{ssec_randst} is
dedicated to the random case with a stable slow manifold and
Subsection~\ref{ssec_bif} to the case of bifurcations.
Sections~\ref{sec_exit} to \ref{sec_bifurcations} contain the proofs of
these results. 

\subsection*{Acknowledgements}
B.\,G. thanks the Forschungsinstitut f\"ur Mathematik at the ETH Z\"urich
and its director Professor Marc Burger for kind hospitality.


\section{Results}   \label{sec_results}


\subsection{Preliminaries}   \label{ssec_prelim}

Let $\cD$ be an open subset of $\R^n\times\R^m$ and $\eps_0>0$
a constant.  We consider slow--fast stochastic differential equations of the
form 
\begin{equation}   \label{rp1}
\begin{split}
\6x_t &= \frac1\eps f(x_t,y_t,\eps)\6t + \frac\sigma{\sqrt\eps}
F(x_t,y_t,\eps)\6W_t, \\
\6y_t &= g(x_t,y_t,\eps)\6t + \sigma' G(x_t,y_t,\eps) \6W_t,
\end{split}
\end{equation}
with drift coefficients $f\in\cC^2(\cD\times[0,\eps_0),\R^n)$ and
$g\in\cC^2(\cD\times[0,\eps_0),\R^m)$, and diffusion coefficients 
$F\in\cC^1(\cD\times[0,\eps_0),\R^{n\times k})$ and  
$G\in\cC^1(\cD\times[0,\eps_0),\R^{m\times k})$.  

We require that $f$, $g$, and all their derivatives up to order $2$ are
uniformly bounded in norm in $\cD\times[0,\eps_0)$, and similarly for $F$,
$G$ and their derivatives. We also assume that $f$ and $g$ satisfy the
usual (local) Lipschitz and bounded-growth conditions which guarantee
existence and pathwise uniqueness of a strong solution
$\{(x_t,y_t)\}_{t\geqs t_0}$ of~\eqref{rp1}. 

The stochastic process $\{W_t\}_{t\geqs 0}$ is a standard $k$-dimensional
Brownian motion on some probability space $(\Omega, \cF, \fP)$, and stochastic
integrals with respect to $\{W_t\}_{t\geqs 0}$ are to be understood as It\^o
integrals. Initial conditions $(x_0,y_0)$ are always assumed to be
square-integrable with respect to $\fP$ and independent of $\{W_t\}_{t\geqs
  0}$. Our assumptions on $f$ and $g$ guarantee the existence of a continuous
version of $\{(x_t,y_t)\}_{t\geqs 0}$. Therefore we may assume that the paths
$\omega\mapsto (x_t(\omega),y_t(\omega))$ are continuous for $\fP$-almost
all $\omega\in\Omega$.  

We introduce the notation $\fP^{\mskip1.5mu t_0,(x_0,y_0)}$ for the law of
the process $\{(x_t,y_t)\}_{t\geqs t_0}$, starting in $(x_0,y_0)$ at time
$t_0$, and use $\E^{\mskip1.5mu t_0,(x_0,y_0)}$ to denote expectations with
respect to $\fP^{\mskip1.5mu t_0,(x_0,y_0)}$. Note that the stochastic
process $\{(x_t,y_t)\}_{t\geqs t_0}$ is a time-homogeneous Markov process. 
Let $\cA\subset\cD$ be Borel-measurable. Assuming $(x_0,y_0)\in\cA$, we
denote by 
\begin{equation}   \label{rp2}
\tau_{\cA}=\inf\bigsetsuch{t\geqs0}{(x_t,y_t)\not\in\cA}
\end{equation}
the first-exit time of $(x_t,y_t)$ from $\cA$. Note that $\tau_{\cA}$ is a
stopping time with respect to the filtration of $(\Omega, \cF, \fP)$
generated by the Brownian motion $\{W_t\}_{t\geqs0}$.

Throughout this work, we use the following notations:

\begin{itemiz}
\item   Let $a$, $b$ be real numbers. We denote by $\intpartplus{a}$,
$a\wedge b$ and $a\vee b$, respectively, the smallest integer greater than
or equal to $a$, the minimum of $a$ and $b$, and the maximum of $a$ and $b$.

\item   By $g(u)=\Order{u}$ we indicate that there exist $\delta>0$ and $K>0$
such that $g(u)\leqs K u$ for all $u\in[0,\delta]$, where $\delta$ and
$K$ of course do not depend on $\eps$, $\sigma$ or $\sigma^\prime$. 

\item   We use $\norm{x}$ to denote the Euclidean norm of $x\in\R^d$ and 
  $\pscal{\cdot}{\cdot}$ for the associated inner product. For a matrix 
$A\in\R^{d_1\times d_2}$, we denote by $\norm{A}$ the corresponding operator
norm. If $A(t)$ is a matrix-valued function defined for $t$ in an interval
$I$, we denote by $\norm{A}_I$ the supremum of $\norm{A(t)}$ over $t\in I$,
and often we write $\norm{A}_\infty$ if the interval is evident from the
context. 

\item   We write $\transpose{A}$ for the transposed of a matrix, and $\Tr A$
for the trace of a square matrix. 

\item   For a given set $B$, we denote by $1_B$ the indicator
function on $B$, defined by $1_B(x) = 1$, if $x\in B$, and
$1_B(x) = 0$, otherwise.

\item   If $\R^n\times\R^m \ni (x,y)\mapsto f(x,y) \in \R^d$ is
differentiable, we write $\partial_x f(x,y)$ and $\partial_y f(x,y)$ to
denote the Jacobian matrices of $x \mapsto f(x,y)$ and $y \mapsto f(x,y)$,
respectively.
\end{itemiz}


\subsection{Deterministic stable case}   \label{ssec_detst}

We start by recalling a few properties of deterministic slow--fast systems of
the form

\begin{equation}   \label{detst1}
\begin{split}
\eps \dot x & = f(x,y,\eps), \\
     \dot y & = g(x,y,\eps). 
\end{split}
\end{equation}

\begin{definition}   \label{def_slowman}
Let $\cD_0 \subset \R^m$ and assume that there exists a (continuous) function
$\xstar : \cD_0 \to \R^n$ such that 
\begin{itemiz}
\item $(\xstar(y),y) \in \cD$ for all $y \in \cD_0$,
\item $f(\xstar(y),y,0)=0$ for all $y \in \cD_0$.
\end{itemiz}
Then the set $\setsuch{(x,y)}{x=\xstar(y), y\in\cD_0}$ is called a\/ {\em slow
  manifold} of the system~\eqref{detst1}. 

Let $\Astar(y) = \partial_x f(\xstar(y),y,0)$. The slow manifold is called
\begin{itemiz}
\item {\em hyperbolic} if all eigenvalues of $\Astar(y)$ have nonzero
  real parts for all $y \in \cD_0$; 
\item {\em uniformly hyperbolic} if all eigenvalues of $\Astar(y)$ have real
  parts uniformly bounded away from zero (for $y \in \cD_0$);
\item {\em asymptotically stable} if all eigenvalues of $\Astar(y)$ have
  negative real parts for all $y \in \cD_0$;
\item {\em uniformly asymptotically stable} if all eigenvalues of
  $\Astar(y)$ have negative real parts, uniformly bounded away from
  zero for $y \in \cD_0$.
\end{itemiz}
\end{definition}
Grad{\v{s}}te{\u\i}n~\cite{Grad} and Tihonov~\cite{Tihonov} have shown that
if $\xstar$ represents a uniformly hyperbolic slow manifold
of~\eqref{detst1}, then the system~\eqref{detst1} admits particular
solutions which remain in a neighbourhood of order~$\eps$ of the slow
manifold. If, moreover, the slow manifold is asymptotically stable, then
the solutions starting in a neighbourhood of order~$1$ of the slow manifold
converge exponentially fast in $t/\eps$ to an $\eps$-neighbourhood of the
slow manifold.

Fenichel~\cite{Fenichel} has given extensions of this result based on a
geometrical approach. If~\eqref{detst1} admits a hyperbolic slow manifold,
then there exists, for sufficiently small $\eps$, an invariant manifold
\begin{equation}   \label{detst2}
x = \bx(y,\eps) = \xstar(y) + \Order{\eps},   \qquad y\in\cD_0.
\end{equation}
Here invariant means that if $y_0\in\cD_0$ and $x_0 = \bx(y_0,\eps)$, then
$x_t=\bx(y_t,\eps)$ as long as $t$ is such that $y_s\in\cD_0$ for all
$s\leqs t$.
We will call the set $\setsuch{(\bx(y,\eps),y)}{y\in\cD_0}$ an {\em
  adiabatic manifold}. It is easy to see from~\eqref{detst1} that
$\bx(y,\eps)$ must satisfy the PDE
\begin{equation}   \label{detst3}
\eps \partial_y \bx(y,\eps) g(\bx(y,\eps),y,\eps) = f(\bx(y,\eps),y,\eps).
\end{equation}
The local existence of the adiabatic manifold follows directly from the
centre manifold theorem. Indeed, we can rewrite System~\eqref{detst1} in the
form 
\begin{equation}   \label{detst4}
\begin{split}
   x^\prime & = f(x,y,\eps), \\
   y^\prime & = \eps g(x,y,\eps), \\
\eps^\prime & = 0,
\end{split}
\end{equation}
where prime denotes derivation with respect to the fast time $t/\eps$. Any
point of the form $(\xstar(y),y,0)$ with $y\in\cD_0$ is an equilibrium
point of~\eqref{detst4}. The linearization of~\eqref{detst4} around such a
point admits $0$ as eigenvalue of multiplicity $m+1$, the $n$ other
eigenvalues being those of $\Astar(y)$, which are bounded away from the
imaginary axis. The centre manifold theorem implies the existence of a
local invariant manifold $x=\bx(y,\eps)$. Fenichel's result shows that this
manifold actually exists for all $y\in\cD_0$.

Being a centre manifold, the adiabatic manifold is not necessarily unique
(though in the present case, $\bx(y,0) = \xstar(y)$ is uniquely defined).
Nevertheless, $\bx(y,\eps)$ has a unique Taylor series in $y$ and $\eps$,
which can be obtained by solving~\eqref{detst3} order by order. The
dynamics on the adiabatic manifold is described by the so-called {\em
reduced equation\/} 
\begin{equation}   \label{detst5}
\dot y = g(\bx(y,\eps),y,\eps) = g(\xstar(y),y,0) + \Order{\eps}.
\end{equation}
If $\xstar(y)$ is uniformly asymptotically stable, $\bx(y,\eps)$ is locally
attractive and thus any solution of~\eqref{detst1} starting sufficiently
close to $\bx(y,\eps)$ converges exponentially fast to a solution
of~\eqref{detst5}.


\subsection{Random stable case}   \label{ssec_randst}

We turn now to the random slow--fast system given by the stochastic
differential equation
\begin{equation}   \label{randst1}
\begin{split}
\6x_t & = \frac1\eps f(x_t,y_t,\eps) \6t + \frac{\sigma}{\sqrt\eps}
F(x_t,y_t,\eps) \6W_t, \\
\6y_t & = g(x_t,y_t,\eps) \6t + \sigma^\prime G(x_t,y_t,\eps) \6W_t,
\end{split}
\end{equation}
where we will assume the following.

\begin{assump}  \label{assump_sde}
For $\sigma=\sigma^\prime=0$, System~\eqref{randst1} admits a
uniformly hyperbolic, asymptotically stable slow manifold
$x=\xstar(y)$, $y\in\cD_0$. 
\end{assump}

By Fenichel's theorem, there exists an adiabatic manifold $x=\bx(y,\eps)$ with
$\bx(y,0) = \xstar(y)$, $y\in\cD_0$. We fix a particular solution
$(\xdet_t,\ydet_t)=(\bx(\ydet_t,\eps),\ydet_t)$ of the deterministic
system. (That is, $\ydet_t$ satisfies the reduced
equation~\eqref{detst5}.) We want to describe the noise-induced deviations
of the sample paths $(x_t,y_t)_{t\geqs 0}$ of~\eqref{randst1} from the
adiabatic manifold. 

It turns out to be convenient to use the transformation
\begin{equation}   \label{randst3}
\begin{split}
x_t & = \bx(\ydet_t+\eta_t,\eps) + \xi_t, \\
y_t & = \ydet_t + \eta_t,
\end{split}
\end{equation}
which yields a system of the form
\begin{equation}   \label{randst4}
\begin{split}
\6\xi_t & = \frac1\eps \hat f(\xi_t,\eta_t,t,\eps) \6t 
+ \frac{\sigma}{\sqrt\eps}
\widehat F(\xi_t,\eta_t,t,\eps) \6W_t, \\
\6\eta_t & = \hat g(\xi_t,\eta_t,t,\eps) \6t + \sigma^\prime \widehat
G(\xi_t,\eta_t,t,\eps) \6W_t, 
\end{split}
\end{equation}
where the new drift and diffusion coefficients are given by
\begin{align}   
\nonumber
\hat f(\xi,\eta,t,\eps) 
& = f(\bx(\ydet_t+\eta,\eps)+\xi,\ydet_t+\eta,\eps) \\
\nonumber
& \hphantom{{}={}}{}- \eps \partial_y \bx(\ydet_t+\eta,\eps)
g(\bx(\ydet_t+\eta,\eps)+\xi,\ydet_t+\eta,\eps)
- \eps\rho^2\sigma^2 r(\xi,\eta,t,\eps), \\
\nonumber
\widehat F(\xi,\eta,t,\eps) 
& = F(\bx(\ydet_t+\eta,\eps)+\xi,\ydet_t+\eta,\eps) \\
\nonumber
&  \hphantom{{}={}}{}- \rho\sqrt\eps \partial_y \bx(\ydet_t+\eta,\eps)
G(\bx(\ydet_t+\eta,\eps)+\xi,\ydet_t+\eta,\eps), \\
\nonumber
\hat g(\xi,\eta,t,\eps) 
& = g(\bx(\ydet_t+\eta,\eps)+\xi,\ydet_t+\eta,\eps)
- g(\bx(\ydet_t,\eps),\ydet_t,\eps), \\
\label{randst5}
\widehat G(\xi,\eta,t,\eps) 
& = G(\bx(\ydet_t+\eta,\eps)+\xi,\ydet_t+\eta,\eps).
\end{align}
Here $r(\xi,\eta,t,\eps)$ stems from the contribution of the diffusion
coefficients in~\eqref{randst1} to the new drift coefficient,
cf.~It\^o's formula. The $l$\/th component of $r(\xi,\eta,t,\eps)$
equals 
\begin{equation}   \label{randst5.5}
\frac12  \Tr\Bigpar{\partial_{yy}\bx_l(\ydet_t,\eps)
G(\bx(\ydet_t+\eta,\eps)+\xi,\ydet_t+\eta,\eps)
\transpose{G(\bx(\ydet_t+\eta,\eps)+\xi,\ydet_t+\eta,\eps)}}, 
\end{equation}
where $\bx_l(y,\eps)$ denotes the $l$\/th component of $\bx(y,\eps)$, and
$\partial_{yy}\bx_l(y,\eps)$ the Hessian matrix of $y\mapsto
\bx_l(y,\eps)$. In the sequel, we will only use the fact that each
component of $r(\xi,\eta,t,\eps)$ is at most of order $m$. 

Note that because of the property~\eqref{detst3} of the adiabatic manifold,
we have $\hat f(0,\eta,t,\eps)=-\eps\rho^2\sigma^2 r(0,\eta,t,\eps)
= \Order{m\eps\rho^2\sigma^2}$. We introduce the notation 
\begin{equation}   \label{randst6}
A(\ydet_t,\eps) 
= \partial_x f(\bx(\ydet_t,\eps),\ydet_t,\eps) 
- \eps \partial_y \bx(\ydet_t,\eps)
\partial_x g(\bx(\ydet_t,\eps),\ydet_t,\eps)
\end{equation}
as an approximation for the linearization of $\hat f$ at $(0,0,t,\eps)$,
where we neglect the contribution of $r(\xi,\eta,t,\eps)$ to the
linearization. Note that for $\eps=0$, we 
have $A(\ydet_t,0) = \partial_x f(\bx(\ydet_t,0),\ydet_t,0) =
\Astar(\ydet_t)$,  so that by Assumption~\ref{assump_sde}, the eigenvalues
of $A(\ydet_t,\eps)$ have negative real parts for sufficiently small
$\eps$.

One of the basic ideas of our approach is to compare the solutions
of~\eqref{randst4} with those of the \lq\lq linear approximation\rq\rq
\begin{equation}   \label{randst7}
\begin{split}
\6\xi^0_t & = \frac1\eps A(\ydet_t,\eps)\xi^0_t \6t + \frac{\sigma}{\sqrt\eps}
F_0(\ydet_t,\eps) \6W_t, \\
\6\ydet_t & = g(\bx(\ydet_t,\eps),\ydet_t,\eps) \6t,
\end{split}
\end{equation}
where $F_0(\ydet_t,\eps) = \widehat F(0,0,t,\eps)$. Note that the definition
of the adiabatic manifold implies $F_0(y,0)=F(\xstar(y),y,0)$. For fixed $t$,
$\xi^0_t$ is a Gaussian random variable with covariance matrix
\begin{equation}   \label{randst8}
\cov(\xi^0_t) = \frac{\sigma^2}{\eps} 
\int_0^t U(t,s) F_0(\ydet_s,\eps) \transpose{F_0(\ydet_s,\eps)}
\transpose{U(t,s)} \6s,  
\end{equation}
where $U(t,s)$ denotes the principal solution of the homogeneous system 
$\eps \dot\xi = A(\ydet_t,\eps)\xi$.

We now observe that $\sigma^{-2} \cov(\xi^0_t)$ is the $X$-variable of a
particular solution of the deterministic slow--fast system
\begin{equation}   \label{randst9}
\begin{split}
\eps \dot X & = A(y,\eps)X + X \transpose{A(y,\eps)} 
+  F_0(y,\eps) \transpose{F_0(y,\eps)}, \\
\dot y & = g(\bx(y,\eps),y,\eps).
\end{split}
\end{equation}
This system admits a slow manifold $X=X^\star(y)$, given by the Lyapunov
equation 
\begin{equation}   \label{randst10}
\Astar(y)X^\star(y) + X^\star(y) \transpose{\Astar(y)} 
+  F_0(y,0) \transpose{F_0(y,0)} =0,
\end{equation}
which is known~\cite{Bellman} to admit the (unique) solution
\begin{equation}   \label{randst11}
X^\star(y) = \int_0^\infty \e^{s\Astar(y)} F_0(y,0) \transpose{F_0(y,0)}
\e^{s\transpose{\Astar(y)}}\6s.
\end{equation}
Moreover, the eigenvalues of the operator $X \mapsto AX+X\transpose{A}$ are
exactly $a_i+a_j$, $1\leqs i,j \leqs n$, where $a_i$ are the eigenvalues of
$A$. Thus the slow manifold $X=X^\star(y)$ is uniformly asymptotically
stable (for small enough $\eps$), so that Fenichel's theorem shows the
existence of an adiabatic manifold
\begin{equation}   \label{randst12}
X = \Xbar(y,\eps) =  X^\star(y) + \Order{\eps}.
\end{equation}
Note that $\Xbar(\ydet_t,\eps)$ is uniquely determined by the \lq\lq
initial\rq\rq\ value $\Xbar(\ydet_0,\eps)$ via the relation
\begin{equation}   \label{randst13}
\Xbar(\ydet_t,\eps) = U(t) \biggbrak{\Xbar(\ydet_0,\eps)
+ \frac{1}{\eps} 
\int_0^t U(s)^{-1} F_0(\ydet_s,\eps) \transpose{F_0(\ydet_s,\eps)}
\transposeinverse{U(s)} \6s} \transpose{U(t)},
\end{equation}
where $U(t)=U(t,0)$ and $\transposeinverse{U(s)}=\transpose{\brak{U(s)^{-1}}}$.

We now introduce the set 
\begin{equation}   \label{randst14}
\cB(h) = 
\bigsetsuch{(x,y)}{y\in\cD_0, \bigpscal{(x-\bx(y,\eps))}{\Xbar(y,\eps)^{-1} 
    (x-\bx(y,\eps))} < h^2},
\end{equation}
assuming that $\Xbar(y,\eps)$ is invertible for all $y\in\cD_0$. The
set $\cB(h)$ is a \lq\lq layer\rq\rq\ around the adiabatic manifold 
$x=\bx(y,\eps)$, with ellipsoidal cross-section determined by
$\Xbar(y,\eps)$. For fixed $t$, the solution $\xi^0_t$ of the linear
approximation~\eqref{randst7} is concentrated (in density) 
in the cross-section of $\cB(\sigma)$ taken at $y_t$. Our first main result
(Theorem~\ref{thm1} below) gives conditions under which the whole sample path
$(x_t,y_t)$ of the original equation~\eqref{randst1} is likely to remain in
such a set $\cB(h)$. By 
\begin{equation}   \label{randst14.1}
\tau_{\cB(h)} = 
\inf\setsuch{t\geqs0}{(x_t,y_t)\not\in\cB(h)}
\end{equation}
we denote the first-exit time of the sample path $(x_t,y_t)$ from
$\cB(h)$. In order to estimate the probability of $\tau_{\cB(h)}$ being small,
we need to assume that $\Xbar(y,\eps)$ and $\Xbar(y,\eps)^{-1}$ are uniformly
bounded in $\cD_0$, which excludes purely multiplicative noise.

\begin{remark}   \label{rem_randst1}
Fix $y$ for the moment. If $X^\star(y)^{-1}$ is bounded, then
$\Xbar(y,\eps)^{-1}$ is bounded for sufficiently small $\eps$. A
sufficient condition for $X^\star(y)^{-1}$ to be bounded is that the symmetric
matrix $F_0(y,0) \transpose{F_0(y,0)}$ be positive definite. This condition
is, however, by no means necessary. In fact, $X^\star(y)$ is singular
if and only if there exists a vector $x\ne0$ such that 
\begin{equation}   \label{randst15}
\transpose{F_0(y,0)} \e^{s\transpose{\Astar(y)}}x = 0
\qquad \forall s\geqs0,
\end{equation}
which occurs if and only if
\begin{equation}   \label{randst16}
\transpose{x} \Astar(y)^l F_0(y,0) =0 
\qquad \forall l=0,1,2,\dots
\end{equation}
Because of the Cayley--Hamilton theorem, this relation holds for all
$l\geqs0$ provided it holds for $l=0,\dots,n-1$. Conversely,
$X^\star(y)$ is {\it non\/}singular if and only if the matrix
\begin{equation}
\label{randst16b}
\bigbrak{F_0(y,0) \quad \Astar(y) F_0(y,0) \quad \dots \quad
\Astar(y)^{n-1} F_0(y,0)}  \in\R^{n\times nk}
\end{equation}
has full rank. This condition on the pair $(\Astar(y),F_0(y,0))$ is known as
\defwd{controllability\/} in control theory, where $X^\star(y)$ is called a
\defwd{controllability Grammian\/}. 

In what follows, we need $(\Astar(y),F_0(y,0))$ to be controllable for all
$y\in\cD_0$, but in addition the smallest eigenvalue of $X^\star(y)$
should be uniformly bounded away from zero. 
\end{remark}

\begin{theorem}   \label{thm1}
Assume that $\norm{\Xbar(y,\eps)}$ and
$\norm{\Xbar(y,\eps)^{-1}}$ are uniformly bounded in $\cD_0$. Choose
a deterministic initial condition $y_0\in\cD_0$, $x_0=\bx(y_0,\eps)$, and let
\begin{equation}   \label{randst17}
\tau_{\cD_0} = \inf\setsuch{s>0}{y_s\not\in\cD_0}.
\end{equation}
Then there exist constants $\eps_0, \Delta_0, h_0>0$ (independent of
the chosen initial condition $y_0$) such that for all $\eps\leqs\eps_0$,
$\Delta\leqs\Delta_0$, $h\leqs h_0$, and all\/ $0<\gamma<1/2$, the following
assertions hold.
\begin{itemize}
\item[(a)] {\em The upper bound:} For all $t>0$,
\begin{equation}   \label{randst18a}
\probin{0,(x_0,y_0)}{\tau_{\cB(h)}< t \wedge \tau_{\cD_0}}
\leqs \cC^+_{n,m,\gamma,\Delta}(t,\eps)
\biggpar{1+\frac{h^2}{\sigma^2}} \e^{-\kappa^ + h^2/\sigma^2}, 
\end{equation}
where 
\begin{equation}    \label{randst18a.5}
\kappa^+ = \gamma \bigbrak{1 - \Order{h} - \Order{\Delta} -
  \Order{m\eps\rho^2} - \bigOrder{\e^{-\text{\it const}/\eps}/(1-2\gamma)}}
\end{equation}
and 
\begin{equation}    \label{randst18a.6}
\cC^+_{n,m,\gamma,\Delta}(t,\eps) = \text{\it const}\
 \frac{(1+t)^2}{\Delta\eps} 
 \Bigbrak{(1-2\gamma)^{-n}  + \e^{n/4} + \e^{m/4}}. 
\end{equation}
\item[(b)] {\em The lower bound:} There exists $t_0>0$ of order $1$ such
for all $t>0$, 
\begin{equation}   \label{randst18b}
\probin{0,(x_0,y_0)}{\tau_{\cB(h)}< t}
\geqs \cC^-_{n,m}(t,\eps,h,\sigma) \e^{-\kappa^- h^2/\sigma^2},
\end{equation}
where 
\begin{equation}    \label{randst18b.5}
\kappa^- = \frac12\bigbrak{1 + \Order{h} 
  + \bigOrder{\e^{-\text{\it const}\,(t\wedge t_0)/\eps}}}
\end{equation}
and 
\begin{equation}    \label{randst18b.6}
\cC^-_{n,m}(t,\eps,h,\sigma) 
= \text{\it const}\ \biggbrak{1 - \biggpar{\e^{n/4} + \frac{\e^{m/4}}{\Delta\eps}}
\e^{-\kappa^- h^2/(2\sigma^2)}}. 
\end{equation}
\item[(c)] {\em General initial conditions:} 
There exist $\delta_0>0$ and a time $t_1$ of order $\eps\abs{\log h}$
such that for all $\delta\leqs \delta_0$, all initial conditions
$(x_0,y_0)$ which satisfy $y_0\in\cD_0$ as well as
$\pscal{\xi_0}{\Xbar(y_0,\eps)^{-1}\xi_0} < \delta^2$, and all $t,t_2$
with $t\geqs t_2\geqs t_1$,
\begin{equation}  \label{randst18c.1}
\biggprobin{0,(\xi_0,0)}{\sup_{t_2 \leqs s \leqs t\wedge\tau_{\cD_0}} 
\bigpscal{\xi_s}{\Xbar(y_s,\eps)^{-1}\xi_s} \geqs h^2}
\leqs \cC^+_{n,m,\gamma,\Delta}(t,\eps)
\biggpar{1+\frac{h^2}{\sigma^2}}  \e^{-\kappa^+ h^2/\sigma^2},
\end{equation}
where $\cC^+_{n,m,\gamma,\Delta}(t,\eps)$ is the same prefactor as
in~\eqref{randst18a}, and 
\begin{equation}  \label{randst18c.2}
\kappa^+ = \gamma \bigbrak{1 - \Order{h} - \Order{\Delta} - \Order{m\eps\rho^2}
- \bigOrder{\delta\e^{-\text{\it const}\,(t_2\wedge1)/\eps}/(1-2\gamma)}}.
\end{equation}
\end{itemize}
Unless explicitly stated, the error terms in the exponents $\kappa^+$ and
$\kappa^-$ are uniform in $t$, but they may depend on the dimensions~$n$
and $m$. 
\end{theorem}

Estimate~\eqref{randst18a} shows that for $h\gg\sigma$, paths starting in
$\cB(h)$ are far more likely to leave this set through the \lq\lq
border\rq\rq\  $\set{y\in\partial\cD_0, \pscal{\xi}{\Xbar(y,\eps)^{-1} \xi}
< h^2}$ than through the \lq\lq sides\rq\rq\  $\set{y\in\interior\cD_0,
\pscal{\xi}{\Xbar(y,\eps)^{-1} \xi} = h^2}$, unless we wait for time spans
exponentially long in $h^2/\sigma^2$. Below we discuss how to
characterize $\tau_{\cD_0}$ more precisely, using information on the
reduced dynamics on the adiabatic manifold. If, for instance, all
deterministic solutions starting in $\cD_0$ remain in this set,
$\tau_{\cD_0}$ will typically be very large. 

The upper bound~\eqref{randst18a} has been designed to yield the best
possible exponent $\kappa^+$, while the prefactor
$\cC^+_{n,m,\gamma,\Delta}$ is certainly not optimal. Note that an estimate
with the same exponent, but with a smaller prefactor holds for the
probability that the endpoint $(x_t,y_t)$ does not lie in $\cB(h)$, cf.\
Corollary~\ref{cor_timetwo}.  The parameters $\Delta$ and $\gamma$ can be
chosen arbitrarily within their intervals of definition. Taking $\Delta$
small and $\gamma$ close to $1/2$ improves the exponent while increasing
the prefactor. A convenient choice is to take $\Delta$ and $1/2-\gamma$ of
order $h$ or $\eps$. The kind of time-dependence of
$\cC^+_{n,m,\gamma,\Delta}$ is probably not optimal, but the fact that
$\cC^+_{n,m,\gamma,\Delta}$ increases with time is to be expected, 
since it reflects the fact that the probability of observing paths making
excursions away from the adiabatic manifold increases with time. As for the
dependence of the prefactor on the dimensions $n$ and $m$, it is due to the
fact that the tails of standard Gaussian random variables show their typical
decay only outside a ball of radius scaling with the square-root of the
dimension. 

The upper bound~\eqref{randst18a} and lower bound~\eqref{randst18b} together
show that the exponential rate of decay of the probability to leave the set
$\cB(h)$ before time $t$ behaves like $h^2/(2\sigma^2)$ in the limit of
$\sigma$, $\eps$ and $h$ going to zero, as one would expect from other
approaches, based for instance on the theory of large deviations. The bounds
hold, however, in a full neighbourhood of $\sigma=\eps=h=0$. 

Finally, Estimate~\eqref{randst18c.1} allows to extend these results to all
initial conditions in a neighbourhood of order $1$ of the adiabatic
manifold. The only difference is that we have to wait for a time of order
$\eps\abs{\log h}$ before the path is likely to have reached the set
$\cB(h)$. After this time, typical paths behave as if they had started on the
adiabatic manifold. 

\begin{remark}   \label{remark_dim}
In Theorem~\ref{thm1}, the error terms in the exponents $\kappa^\pm$ grow 
with the norms $\norm{f}$ and $\norm{g}$, and thus depend in general 
on the dimensions~$n$ and~$m$. If the SDE~\eqref{randst1} describes a large 
number of coupled similar subsystems (e.\,g.\ coupled oscillators), the error 
terms will {\em not\/} depend on the number of subsystems if, for instance, 
each one is coupled only to a finite number of neighbours. In mean-field type 
models, the error terms will be bounded if the interaction is properly scaled 
with the number of subsystems. 
\end{remark}

The behaviour of typical paths depends essentially on the dynamics of the
reduced deterministic system~\eqref{detst5}. In fact, in the proof of
Theorem~\ref{thm1}, we use the fact that $y_t$ does not differ too much
from $\ydet_t$ on timescales of order~$1$ (see Lemma~\ref{l_eta}). There
are thus two main possibilities to be considered:

\begin{itemiz}
\item   either the reduced flow is such that $\ydet_t$ reaches the boundary
of $\cD_0$ in a time of order~$1$ (for instance, $\ydet_t$ may approach a
bifurcation set of the slow manifold); then $y_t$ is likely to leave
$\cD_0$ as well;

\item   or the reduced flow is such that $\ydet_t$ remains in $\cD_0$ for
all times $t\geqs0$; in that case, paths can only leave $\cB(h)$ due to the
influence of noise, which we expect to be unlikely on subexponential
timescales.
\end{itemiz}

We will discuss the first situation in more detail in
Subsection~\ref{ssec_bif}. In both situations, it is desirable to have
a more precise description of the deviation $\eta_t$ of the slow variable
$y_t$ from its deterministic counterpart $\ydet_t$, in order to achieve a
better control of the first-exit time $\tau_{\cD_0}$. 

The following coupled system gives a better approximation of the dynamics
of~\eqref{randst4} than the system~\eqref{randst7}:
\begin{equation}   \label{randst21}
\begin{split}
\6\xi^0_t & = \frac1\eps A(\ydet_t,\eps)\xi^0_t \6t + \frac{\sigma}{\sqrt\eps}
F_0(\ydet_t,\eps) \6W_t, \\
\6\eta^0_t & = \bigbrak{B(\ydet_t,\eps)\eta^0_t + C(\ydet_t,\eps)\xi^0_t} \6t 
+ \sigma^\prime G_0(\ydet_t,\eps) \6W_t,
\end{split}
\end{equation}
where $G_0(\ydet_t,\eps) = \widehat G(0,0,t,\eps) =
G(\bx(\ydet_t,\eps),\ydet_t,\eps)$ and the Jacobian matrices $B$ and $C$ are
given by
\begin{align}   
\nonumber
B(\ydet_t,\eps) &= \sdpar{\hat g}\eta(0,0,t,\eps) \\
\label{randst22}
&= C(\ydet_t,\eps) \sdpar\bx y(\ydet_t,\eps)
+ \sdpar gy(\bx(\ydet_t,\eps),\ydet_t,\eps), \\
\nonumber
C(\ydet_t,\eps) &= \sdpar{\hat g}\xi(0,0,t,\eps) \\
\label{randst23}
&= \sdpar gx(\bx(\ydet_t,\eps),\ydet_t,\eps). 
\end{align}
The coupled system~\eqref{randst21} can be written in compact form as
\begin{equation}   \label{randst24}
\6\z^0_t = \cA(\ydet_t,\eps)\z^0_t \6t + \sigma\cF_0(\ydet_t,\eps) \6W_t,
\end{equation} 
where $\transpose{(\z^0)} = (\transpose{(\xi^0)},\transpose{(\eta^0)})$ and
\begin{equation}   \label{randst25}
\cA(\ydet_t,\eps) = 
\begin{pmatrix}
\lower10pt\vbox{\kern0pt} \frac1\eps A(\ydet_t,\eps) & 0 \\
\raise10pt\vbox{\kern0pt} C(\ydet_t,\eps) & B(\ydet_t,\eps)
\end{pmatrix},
\qquad
\cF_0(\ydet_t,\eps) = 
\begin{pmatrix}
\lower10pt\vbox{\kern0pt} \frac1{\sqrt\eps} F_0(\ydet_t,\eps) \\
\raise10pt\vbox{\kern0pt} \rho G_0(\ydet_t,\eps)
\end{pmatrix}.
\end{equation}
The solution of the linear SDE~\eqref{randst24} is given by 
\begin{equation}   \label{randst26}
\z^0_t = \cU(t) \z_0 + \sigma\int_0^t \cU(t,s) \cF_0(\ydet_s,\eps) \6W_s,
\end{equation}
where $\cU(t,s)$ denotes the principal solution of the homogeneous system
$\dot\z = \cA(\ydet_t,\eps)\z$. It can be written in the form 
\begin{equation}   \label{randst27}
\cU(t,s) = 
\begin{pmatrix}
U(t,s) & 0 \\ S(t,s) & V(t,s)
\end{pmatrix},
\end{equation}
where $U(t,s)$ and $V(t,s)$ denote, respectively, the fundamental solutions
of $\eps\dot\xi = A(\ydet_t,\eps)\xi$ and $\dot\eta = B(\ydet_t,\eps)\eta$,
while 
\begin{equation}   \label{randst28}
S(t,s) = \int_s^t V(t,u) C(\ydet_u,\eps) U(u,s) \6u.
\end{equation}
The Gaussian process $\z^0_t$ has a covariance matrix of the form 
\begin{align}
\nonumber
\cov(\z^0_t) &= \sigma^2
\int_0^t \cU(t,s) \cF_0(\ydet_s,\eps) \transpose{\cF_0(\ydet_s,\eps)}
\transpose{\cU(t,s)} \6s \\  
\label{randst29}
&= \sigma^2 
\begin{pmatrix}
X(t) & Z(t) \\ \transpose{Z(t)} & Y(t)
\end{pmatrix}.
\end{align}
The matrices $X(t)\in\R^{n\times n}$, $Y(t)\in\R^{m\times m}$ and
$Z(t)\in\R^{n\times m}$ are a particular solution of the
following slow--fast system, which generalizes~\eqref{randst9}:
\begin{equation}   \label{randst30}
\begin{split}
\eps \dot X & = A(y,\eps)X + X \transpose{A(y,\eps)} 
+  F_0(y,\eps) \transpose{F_0(y,\eps)}, \\
\eps \dot Z & = A(y,\eps)Z + \eps Z \transpose{B(y,\eps)} 
+ \eps X \transpose{C(y,\eps)} 
+ \sqrt\eps\rho F_0(y,\eps) \transpose{G_0(y,\eps)}, \\
\dot Y & = B(y,\eps)Y + Y \transpose{B(y,\eps)} 
+ C(y,\eps)Z + \transpose Z \transpose{C(y,\eps)} 
+ \rho^2 G_0(y,\eps) \transpose{G_0(y,\eps)}, \\
\dot y & = g(\bx(y,\eps),y,\eps).
\end{split}
\end{equation}
This system admits a slow manifold given by 
\begin{align}
\nonumber
X &= X^\star(y), \\
\label{randst31}
Z &= Z^\star(y,\eps) 
= -\sqrt\eps \rho A(y,\eps)^{-1} F_0(y,\eps) \transpose{G_0(y,\eps)} 
+ \Order{\eps},
\end{align}
where $X^\star(y)$ is given by~\eqref{randst11}. It is straightforward to
check that this manifold is uniformly asymptotically stable for sufficiently
small $\eps$, so that Fenichel's theorem yields the existence of an
adiabatic manifold $X = \Xbar(y,\eps)$, $Z = \Zbar(y,\eps)$, at a distance
of order~$\eps$ from the slow manifold. This manifold attracts nearby
solutions of~\eqref{randst30} exponentially fast, and thus asymptotically,
the expectations of $\xi^0_t \transpose{(\xi^0_t)}$ and $\xi^0_t
\transpose{(\eta^0_t)}$ will be close, respectively, to
$\sigma^2\Xbar(\ydet_t,\eps)$ and $\sigma^2\Zbar(\ydet_t,\eps)$. 

In general, the matrix $Y(t)$ cannot be expected to approach some
asymptotic value depending only on $\ydet_t$ and $\eps$. In fact, if the
deterministic orbit $\ydet_t$ is repelling, $\norm{Y(t)}$ can grow
exponentially fast. In order to measure this growth, we introduce the
functions  
\begin{align}   
\label{randst34a}
\chi^{(1)}(t) & = \sup_{0 \leqs s\leqs t} \int_0^s \Bigpar{\sup_{u\leqs v\leqs
    s} \norm{V(s,v)}} \6u,  \\
\label{randst34b}
\chi^{(2)}(t) & = \sup_{0 \leqs s\leqs t} \int_0^s \Bigpar{\sup_{u\leqs v\leqs
    s} \norm{V(s,v)}^2} \6u.
\end{align}
The solution of~\eqref{randst30} with initial condition $Y(0)=Y_0$
satisfies 
\begin{align}   
\label{randst32}
Y(t;Y_0) ={}& V(t)Y_0 \transpose{V(t)} \\
\nonumber
&{}+ \rho^2 \int_0^t V(t,s)  G_0(\ydet_s,\eps)
\transpose{G_0(\ydet_s,\eps)} \transpose{V(t,s)} \6s +
\Order{(\eps+\rho\sqrt\eps)\chi^{(2)}(t)}. 
\end{align}
We thus define an \lq\lq asymptotic\rq\rq\  covariance matrix
$\cZbar(t) = \cZbar(t;Y_0,\eps)$ by
\begin{equation}   \label{randst33}
\cZbar(t;Y_0,\eps) = 
\begin{pmatrix}
\lower8pt\vbox{\kern0pt}
\Xbar(\ydet_t,\eps)\phantom{^T} & \Zbar(\ydet_t,\eps) \\
\raise8pt\vbox{\kern0pt}
\transpose{\Zbar(\ydet_t,\eps)} & Y(t;Y_0)
\end{pmatrix},
\end{equation}
and use $\cZbar(t)^{-1}$ to characterize the ellipsoidal region in which
$\z(t)$ is concentrated.

\begin{theorem}
\label{thm2}
Assume that $\norm{\Xbar(\ydet_s,\eps)}$ and
$\norm{\Xbar(\ydet_s,\eps)^{-1})}$ are uniformly bounded for $0\leqs s\leqs t$
and that $Y_0$ has been chosen in such a way that 
$\norm{Y(s)^{-1}} = \Order{1/(\rho^2+\eps)}$ for $0\leqs s\leqs t$. Fix an
initial condition $(x_0,y_0)$ with $y_0\in\cD_0$ and 
$x_0=\bx(y_0,\eps)$, and let $t$ be such that $\ydet_s\in\cD_0$ for all
$s\leqs t$. Define 
\begin{equation}
\label{randst35a}
R(t) = \norm{\cZbar}_{[0,t]}
\Bigbrak{1+\Bigpar{1+\norm{Y^{-1}}_{[0,t]}^{1/2}} \chi^{(1)}(t) +
  \chi^{(2)}(t)}.
\end{equation}
There exist constants $\eps_0, \Delta_0, h_0>0$, independent of $Y_0$, $y_0$
and $t$, such that 
\begin{equation}   \label{randst36}
\Bigprobin{0,(0,0)}{\sup_{0\leqs s\leqs t\wedge\tau_{\cD_0}} 
\bigpscal{\z_u}{\cZbar(u)^{-1}\z_u}\geqs h^2} 
\leqs \cC_{n+m,\gamma,\Delta}(t,\eps) \e^{-\kappa h^2/\sigma^2}
\end{equation}
holds, whenever $\eps\leqs \eps_0$, $\Delta\leqs\Delta_0$, $h\leqs
h_0 R(t)^{-1}$ and $0<\gamma<1/2$. Here
\begin{align}   \label{randst37a}
\cC_{n+m,\gamma,\Delta}(t,\eps) &= \text{\it const\;}
\biggintpartplus{\frac{t}{\Delta\eps}}
\biggbrak{\biggpar{\frac{1}{1-2\gamma}}^{(n+m)/2} + \e^{(n+m)/4}}, \\
\kappa & = \vphantom{\biggpar{\frac1{1-2\gamma}}^{(n+m)/2}}
\gamma \Bigbrak{1 - \bigOrder{\eps + \Delta + h R(t)}}. 
\end{align}
\end{theorem}

Let us first consider timescales of order $1$. Then the functions
$\norm{\cZbar}_{[0,t]}$, $\chi^{(1)}(t)$ and $\chi^{(2)}(t)$ are at most
of order $1$, and $\norm{Y(t)^{-1}}$ remains of the same order as
$\norm{Y_0^{-1}}$. The probability~\eqref{randst36} becomes small as soon as
$h\gg\sigma$. Because of the restriction $h\leqs h_0 R(t)^{-1}$, the result
is useful provided $\norm{Y^{-1}}_{[0,t]} \ll \sigma^{-2}$. In order to obtain
the optimal concentration result, we have to choose $Y_0$ according to two
opposed criteria. On the one hand, we would like to choose $Y_0$ as small as
possible, so that the set $\bigpscal{\z_u}{\cZbar(u)^{-1}\z_u}< h^2$ is
small. On the other hand, $\norm{Y_0^{-1}}$ must not exceed certain bounds for
Theorem~\ref{thm2} to be valid. Thus we require that
\begin{equation}   \label{randst41}
Y_0 > \bigbrak{\sigma^2 \vee (\rho^2+\eps)} \one_m.
\end{equation}
Because of the Gaussian decay of the probability~\eqref{randst36} in
$\sigma/h$, we can interpret the theorem by saying that the typical
spreading of paths in the $y$-direction is of order
$\sigma(\rho+\sqrt\eps)$ if $\sigma<\rho+\sqrt\eps$ and of order $\sigma^2$
if $\sigma>\rho+\sqrt\eps$.

The term $\rho$ is clearly due to the intensity $\sigma^\prime = \rho\sigma$
of the noise acting on the slow variable. It prevails if
$\rho>\sigma\vee\sqrt\eps$. The term $\sqrt\eps$ is due to the linear part of
the coupling between slow and fast variables, while the behaviour in $\sigma^2$
observed when $\sigma>\rho+\sqrt\eps$ can be traced back to the {\em
nonlinear} coupling between slow and fast variables. 

For longer timescales, the condition $h\leqs h_0 R(t)^{-1}$ obliges us to
take a larger $Y_0$, while $Y(t)$ typically grows with time. If the largest
Lyapunov exponent of the deterministic orbit $\ydet_t$ is positive,
this growth is exponential in time, so that the spreading of paths
along the adiabatic manifold will reach order $1$ in a time of order
$\log\abs{\sigma\vee(\rho^2+\eps)}$. 

\begin{remark}   \label{rem_reduced}
Consider the \defwd{reduced stochastic system}
\begin{equation}   \label{randst41a}
\6y^0_t = g(\bx(y^0_t,\eps),y^0_t,\eps) \6t 
+ \sigma^\prime G(\bx(y^0_t,\eps),y^0_t,\eps) \6W_t
\end{equation}
obtained by setting $x$ equal to $\bx(y,\eps)$ in~\eqref{randst1}. One may
wonder whether $y^0_t$ gives a better approximation of 
$y_t$ than $\ydet_t$ in the case $\sigma'>0$. In fact, one can show that 
\begin{align}   
\nonumber
\Bigprobin{0,(0,0)}{\sup_{0\leqs s\leqs t\wedge\tau_{\cD_0}}
\bignorm{y^0_s-\ydet_s}\geqs h_1}
&\leqs c\Bigpar{1+t}\e^{m/4} \exp\Bigset{-\frac{\kappa_1h_1^2}
  {(\sigma^\prime)^2(1+\chi^{(2)}(t))}}, \\ 
\label{randst41b}
\Bigprobin{0,(0,0)}{\sup_{0\leqs s\leqs t\wedge\tau_{\cB(h)}}
\bignorm{y_s-y^0_s}\geqs h} 
&\leqs  c\Bigpar{1+t}\e^{m/4} \exp\Bigset{-\frac{\kappa_1h_1^2}
  {(\sigma^\prime)^2(1+\chi^{(2)}(t))}} \\
\nonumber
& \;{}+ c\Bigpar{1{}+\frac t\eps} \e^{m/4} \exp\Bigset{-\frac{\kappa_2h^2}
{\brak{(\sigma^\prime)^2h^2+\sigma^2\eps}(1+\chi^{(2)}(t))}}
\end{align}
holds for all $h,h_1$ up to order $\chi^{(1)}(t)^{-1}$ and some positive
constants $c,\kappa_1,\kappa_2$. (The proofs can be adapted from the proof of
Lemma~\ref{l_eta}). This shows that the typical spreading of $y^0_t$ around
$\ydet_t$ is of order $\sigma^\prime(1+\chi^{(2)}(t)^{1/2}) =
\rho\sigma(1+\chi^{(2)}(t)^{1/2})$, while the typical deviation of paths
$y^0_t$ of the reduced system from paths $y_t$ of the original system is of
order $\sigma\sqrt\eps(1+\chi^{(2)}(t)^{1/2})$. Thus for $\rho>\sqrt\eps$, the
reduced stochastic system gives a better approximation of the dynamics than
the deterministic one. 
\end{remark}

If $V(t)$ has no eigenvalues outside the unit circle, the spreading of paths
will grow more slowly. As an important particular case, let us consider the
situation where $\ydet_t$ is an asymptotically stable periodic orbit with
period $T$, entirely contained in $\cD_0$ (and not too close to its
boundary). Then all coefficients in~\eqref{randst21} depend periodically on
time, and, in particular, Floquet's theorem allows us to write
\begin{equation}   \label{randst42}
V(t) = P(t)\e^{\Lambda t},
\end{equation}
where $P(t)$ is a $T$-periodic matrix. The asymptotic stability of the orbit
means that all eigenvalues but one of the monodromy matrix $\Lambda$ have
strictly negative real parts, the last eigenvalue, which corresponds to
translations along the orbit, being $0$. In that case, $\chi^{(1)}(t)$ and
$\chi^{(2)}(t)$ grow only linearly with time, so that the spreading of paths
in the $y$-direction remains small on timescales of order
$1/(\sigma\vee(\rho^2+\eps))$. 

In fact, we even expect this spreading to occur mainly along the periodic
orbit, while the paths remain confined to a neighbourhood of the orbit on
subexponential timescales. To see that this is true, we can use a new set
of variables in the neighbourhood of the orbit. In order not to introduce too
many new notations, we will replace $y$ by $(y,z)$, where $y\in\R^{m-1}$
describes the degrees of freedom transversal to the orbit, and $z\in\R$
parametrizes the motion along the orbit. In fact, we can use an equal-time
parametrization of the orbit, so that $\dot z=1$ on the orbit, 
i.\,e., we have $\zdet_t = t \pmod{T}$. 
The SDE takes the form 
\begin{equation}   \label{randst43}
\begin{split}
\6x_t & = \frac1\eps f(x_t,y_t,z_t,\eps) \6t + \frac{\sigma}{\sqrt\eps}
F(x_t,y_t,z_t,\eps) \6W_t, \\
\6y_t & = g(x_t,y_t,z_t,\eps) \6t + \sigma^\prime G(x_t,y_t,z_t,\eps) \6W_t,
\vphantom{ \frac{\sigma}{\sqrt\eps}}
\\
\6z_t & = \bigbrak{1+h(x_t,y_t,z_t,\eps)}\6t 
+ \sigma^\prime H(x_t,y_t,z_t,\eps) \6W_t,
\vphantom{ \frac{\sigma}{\sqrt\eps}}
\end{split}
\end{equation}
where $h = \Order{\norm{y_t}^2 + \norm{x_t-\xdet_t}^2}$ and the Floquet
multipliers associated with the periodic matrix  $\sdpar
gy(\xdet_t,0,\zdet_t,\eps)$ are strictly smaller than one in modulus. 
As linear approximation of the dynamics of
$(\xi_t,\eta_t)=(x_t-\xdet_t,y_t-\ydet_t)=(x_t-\xdet_t,y_t)$ we take  
\begin{equation}   \label{randst44}
\begin{split}
\6\xi^0_t & = \frac1\eps A(\zdet_t,\eps)\xi^0_t \6t + \frac{\sigma}{\sqrt\eps}
F_0(\zdet_t,\eps) \6W_t, \\
\vphantom{ \frac{\sigma}{\sqrt\eps}}
\6\eta^0_t & = \bigbrak{B(\zdet_t,\eps)\eta^0_t + C(\zdet_t,\eps)\xi^0_t} \6t 
+ \sigma^\prime G_0(\zdet_t,\eps) \6W_t, \\
\vphantom{ \frac{\sigma}{\sqrt\eps}}
\6z^0_t & = \6t + \sigma^\prime H_0(\zdet_t,\eps) \6W_t,
\end{split}
\end{equation}
which depends periodically on time. One can again compute the covariance
matrix of the Gaussian process $(\xi^0_t,\eta^0_t,z^0_t)$ as a function of
the principal solutions $U$ and $V$ associated with $A$ and $B$. In
particular, the covariance matrix $Y(t)$ of $\eta^0_t$ still obeys the ODE 
\begin{equation}   \label{randst45}
\dot Y  = B(\zdet,\eps)Y + Y \transpose{B(\zdet,\eps)} 
+ C(\zdet,\eps)\Zbar + \transpose \Zbar \transpose{C(\zdet,\eps)} 
+ \rho^2 G_0(\zdet,\eps) \transpose{G_0(\zdet,\eps)}. 
\end{equation}
This is now a linear, inhomogeneous ODE with time-periodic coefficients. It
is well known that such a system admits a unique periodic solution
$\Yper_t$, which is of order $\rho^2+\eps$ since $\Zbar$ is of order
$\rho\sqrt\eps+ \eps$ and $\rho^2G_0\transpose{G_0}$ is of order $\rho^2$. We
can thus define an asymptotic covariance matrix $\cZbar(t)$ of
$(\x^0_t,\y^0_t)$, which depends periodically on time. If
$\z_t=(\x_t,\y_t)$, Theorem~\ref{thm2} shows that
on timescales of order $1$ (at least), the paths $\z_t$ are concentrated in
a set of the form $\pscal{\z_t}{\cZbar(t)^{-1}\z_t}<h^2$, while $z_t$
remains $h$-close to $\zdet_t$. 

On longer timescales, the distribution of paths will be smeared out along
the periodic orbit. However, the same line of reasoning as in
Section~\ref{ssec_timetwo}, based on a comparison with different
deterministic solutions on successive time intervals of order $1$, can be
used to show that $\z_t$ remains concentrated in the set
$\pscal{\z_t}{\cZbar(t)^{-1}\z_t}<h^2$ up to exponentially long
timescales.


\subsection{Bifurcations}   \label{ssec_bif}

In the previous section, we have assumed that the slow manifold
$x=x^\star(y)$ is uniformly asymptotically stable for $y\in\cD_0$. We
consider now the situation arising when the reduced deterministic flow 
causes $\ydet_t$ to leave $\cD_0$, and to approach a bifurcation point of
the slow manifold.

We call $(\hat x,\hat y)$ a bifurcation point of the deterministic system 
\begin{equation}
\label{bif1}
\begin{split}
\eps\dot{x} &= f(x,y,\eps), \\
\dot{y} &= g(x,y,\eps), 
\end{split}
\end{equation}
if $f(\hat x,\hat y,0)=0$ and $\sdpar fx(\hat x,\hat y,0)$ has $q$
eigenvalues on the imaginary axis, $q\in\set{1,\dots,n}$. We consider here
the situation where $q<n$ and the other $n-q$ eigenvalues have strictly
negative real parts. 

The most generic cases are the saddle--node bifurcation (where $q=1$),
corresponding to a fold in the slow manifold, and the Hopf bifurcation
(where $q=2$), in which the slow manifold changes stability, while
absorbing or expelling a family of periodic orbits. In these two cases, the
set of bifurcation values $\hat y$ typically forms a codimension-1
submanifold of $\R^m$. 

The dynamics of the deterministic slow--fast system~\eqref{bif1} in a
neighbourhood of the bifurcation point $(\hat x,\hat y)$ can again be analyzed
by a centre-manifold reduction. Introduce coordinates $(x^-,z)$ in $\R^n$, with
$x^-\in\R^{n-q}$ and $z\in\R^q$, in which the matrix $\sdpar fx(\hat x,\hat y,
0)$ becomes block-diagonal, with a block $A^-\in\R^{(n-q)\times(n-q)}$ having
eigenvalues in the left half-plane, and a block $A^0\in\R^{q\times q}$ having
eigenvalues on the imaginary axis. On the fast timescale
$t/\eps$,~\eqref{bif1} can be rewritten as  
\begin{equation}   \label{bif2}
\begin{split}
(x^-)^\prime & = f^-(x^-,z,y,\eps), \\
    z^\prime & = f^0(x^-,z,y,\eps), \\
    y^\prime & = \eps g(x^-,z,y,\eps), \\
 \eps^\prime & = 0,
\end{split}
\end{equation}
which admits $(\hat x^-,\hat z,\hat y,0)$ as an equilibrium point. The
linearization at this point has $q+m+1$ eigenvalues on the imaginary axis
(counting multiplicity), which correspond to the directions $z, y$ and
$\eps$. In other words, $z$ has become a slow variable near the bifurcation
point. 

The centre manifold theorem implies the existence, for sufficiently small
$\eps$ and $(z,y)$ in a neighbourhood $\cN$ of $(\hat z,\hat y)$, of a locally
attracting invariant manifold $x^-=\bxm(z,y,\eps)$, with $\bxm(\hat
z,\hat y,0) = \hat x$. $\bxm$ plays the same r\^ole the adiabatic manifold
played in the stable case, and the dynamics on $\bxm$ is governed by the
reduced equation 
\begin{equation}
\label{bif3}
\begin{split}
\eps\dot{z} &= f^0(\bxm(z,y,\eps),z,y,\eps), \\
\dot{y} &= g(\bxm(z,y,\eps),z,y,\eps).  
\end{split}
\end{equation}
The function $\bxm(z,y,\eps)$ solves the PDE
\begin{align}
\nonumber
f^-(\bxm(z,y,\eps),z,y,\eps) 
={}& \sdpar{\bxm}z(z,y,\eps) f^0(\bxm(z,y,\eps),z,y,\eps) \\
&{}+ \eps \sdpar{\bxm}y(z,y,\eps) g(\bxm(z,y,\eps),z,y,\eps). 
\label{bif4}
\end{align}

Let us now turn to random perturbations of the slow--fast
system~\eqref{bif1}. In the variables $(x^-,z,y)$, the perturbed system can be
written as  
\begin{equation}   \label{bif5}
\begin{split}
\6{x^-}_t & = \frac1\eps f^-(x^-_t,z_t,y_t,\eps) \6t + \frac{\sigma}{\sqrt\eps}
F^-(x^-_t,z_t,y_t,\eps) \6W_t, \\
\6z_t & = \frac1\eps f^0(x^-_t,z_t,y_t,\eps) \6t + \frac{\sigma}{\sqrt\eps}
F^0(x^-_t,z_t,y_t,\eps) \6W_t, \\
\6y_t & = 
\vphantom{\frac{1}{\eps}}
g(x^-_t,z_t,y_t,\eps) \6t + \sigma^\prime G(x^-_t,z_t,y_t,\eps)
\6W_t.
\end{split}
\end{equation}
The noise-induced deviation of $x^-_t$ from the adiabatic manifold is
described by the variable $\xi^-_t=x^-_t - \bxm(z_t,y_t,\eps)$, which
obeys an SDE of the form 
\begin{equation}
\label{bif6}
\6\xi^-_t = \frac1\eps \hat f^-(\xi^-_t,z_t,y_t,t,\eps) \6t +
\frac{\sigma}{\sqrt\eps} 
\widehat F^-(\xi^-_t,z_t,y_t,\eps) \6W_t, 
\end{equation}
with, in particular, 
\begin{align}   
\nonumber
\hat f^-(\xi^-,z,y,t,\eps) 
& = f^-(\bxm(z,y,\eps)+\xi^-,z,y,\eps) 
- \partial_z \bxm(z,y,\eps)
f^0(\bxm(z,y,\eps)+\xi^-,z,y,\eps) \\
\label{bif7}
& \hphantom{{}={}}{}- \eps \partial_y \bxm(z,y,\eps)
g(\bxm(z,y,\eps)+\xi^-,z,y,\eps)
-\sigma^2 r^-(\xi^-,z,y,t,\eps),  
\end{align}
where $r^-(\xi^-,z,y,t,\eps)$ is at most of order $m+q$. Note
that~\eqref{bif4} implies that 
\begin{equation}   \label{bif7.5}
\hat f^-(0,z,y,t,\eps) = -\sigma^2  r^-(0,z,y,t,\eps) =
\Order{(m+q)\sigma^2}.
\end{equation} 
We further define the matrix 
\begin{align}   
\nonumber
A^-(z,y,\eps) 
& = \partial_x f^-(\bxm(z,y,\eps),z,y,\eps) 
- \partial_z \bxm(z,y,\eps)
\partial_x f^0(\bxm(z,y,\eps),z,y,\eps) \\
\label{bif8}
& \hphantom{{}={}}{}- \eps \partial_y \bxm(z,y,\eps)
\partial_x g(\bxm(z,y,\eps),z,y,\eps)
\end{align}
as an approximation to $\partial_\xi\hat f^-(0,z,y,t,\eps)$, where we
neglect the contribution of $r^-$.
Since $A^-(\hat z,\hat y,0) = A^-$, the eigenvalues of $A^-(z,y,\eps)$ have
uniformly negative real parts, provided we take the neighbourhood $\cN$ and
$\eps$ small enough. 

Consider now the \lq\lq linear approximation\rq\rq\ 
\begin{equation}   \label{bif9}
\begin{split}
\6\xi^0_t & = \frac1\eps A^-(\zdet_t,\ydet_t,\eps)\xi^0_t \6t 
+ \frac{\sigma}{\sqrt\eps} F_0^-(\zdet_t,\ydet_t,\eps) \6W_t, \\
\6\zdet_t & = \frac1\eps f^0(\bxm(\zdet,\ydet,\eps),\zdet_t,\ydet_t,\eps) \6t,
\\ 
\6\ydet_t & = 
\vphantom{\frac1\eps}
g(\bxm(\zdet,\ydet,\eps),\zdet_t,\ydet_t,\eps) \6t
\end{split}
\end{equation}
of~\eqref{bif5}--\eqref{bif6}, where 
$F_0^-(z,y,\eps)=\widehat F^-(0,z,y,\eps)$. Its solution $\xi^0_t$ has a
Gaussian distribution with covariance matrix
\begin{equation}
\label{bif10}
\cov(\xi^0_t) = \frac{\sigma^2}{\eps} 
\int_0^t U^-(t,s) F_0^-(\zdet_s,\ydet_s,\eps) 
\transpose{F_0^-(\zdet_s,\ydet_s,\eps)} \transpose{U^-(t,s)} \6s,  
\end{equation}
where $U^-$ is the fundamental solution of $\eps\dot\xi^0=A^-\xi^0$. Note that
$\sigma^{-2}\cov(\xi^0_t)$ is the $X^-$-variable of a particular solution of
the slow--fast system 
\begin{equation}   \label{bif11}
\begin{split}
\eps \dot X^- & = A^-(z,y,\eps)X^- + X^- \transpose{A^-(z,y,\eps)} 
+  F_0^-(z,y,\eps) \transpose{F_0^-(z,y,\eps)}, \\
\eps\dot z & = f^0(\bxm(z,y,\eps),z,y,\eps), \\
\dot y & = g(\bxm(z,y,\eps),z,y,\eps),
\end{split}
\end{equation}
which admits an invariant manifold $X^- = \Xbarm(z,y,\eps)$ for
$(z,y)\in\cN$. We thus expect the paths to be concentrated in a set 
\begin{equation}
\label{bif12}
\cB^-(h) = \bigsetsuch{(x^-,z,y)}
{(z,y)\in\cN,
  \bigpscal{x^--\bxm(z,y,\eps)}{\Xbarm(z,y,\eps)^{-1}(x^--\bxm(z,y,\eps))} 
< h^2}. 
\end{equation}
The following theorem shows that this is indeed the case, as long as
$(z_t,y_t)$ remains in $\cN$. 

\begin{theorem}
\label{thm3}
Assume that $\norm{\Xbarm(z,y,\eps)}$ and $\norm{\Xbarm(z,y,\eps)^{-1})}$
are uniformly bounded in $\cN$. Choose a deterministic initial condition
$(z_0,y_0)\in\cN$, $x_0^-=\bxm(z_0,y_0,\eps)$, and let
\begin{equation}   \label{bif13}
\tau_{\cN} = \inf\setsuch{s>0}{(z_s,y_s)\not\in\cN}.
\end{equation}
Then there exist constants $h_0>0$, $\Delta_0>0$ and $\nu\in(0,1]$ such that
for all $h \leqs h_0$, all $\Delta\leqs\Delta_0$ and all $0<\gamma<1/2$,
\begin{equation}   \label{bif14}
\probin{0,(x^-_0,z_0,y_0)}{\tau_{\cB^-(h)}< t \wedge \tau_{\cN}}
\leqs \cC_{n,m,q,\gamma,\Delta}(t,\eps) \biggpar{1+\frac{h^2}{\sigma^2}}
\e^{-\kappa h^2/\sigma^2}, 
\end{equation}
provided $\eps\abs{\log (h(1-2\gamma))}\leqs 1$. Here 
\begin{align}
\kappa &= \gamma\bigbrak{1 - \Order{\Delta} - \Order{h^\nu (1-2\gamma)^{1-\nu}
\abs{\log (h(1-2\gamma))}}},
\label{bif16}
\\
\label{bif15} 
\cC_{n,m,q,\gamma,\Delta}(t,\eps) &
= \text{\it const}\, \Bigpar{1 +\frac{t}{\Delta\eps}}\Bigpar{1 + \frac t\eps}
\Bigbrak{(1-2\gamma)^{-(n-q)}  + \e^{(n-q)/4} + \e^{m/4} + \e^{q/4}}.
\end{align}
\end{theorem}

The exponent $\nu$ is related to the maximal rate of divergence of solutions
of the reduced system~\eqref{bif3}, see Subsection~\ref{ssec_bifexit}.

This result shows that on timescales of order $1$ (and larger if,
e.\,g., $\cN$ is positively invariant), paths are likely to remain in a small
neighbourhood of the adiabatic manifold $x^-=\bxm(z,y,\eps)$. 
The dynamics will thus
be essentially governed by the behaviour of the \lq\lq slow\rq\rq\ variables
$z$ and $y$. 

In fact, it seems plausible that the dynamics of~\eqref{bif5} will be well
approximated by the dynamics of the \defwd{reduced stochastic system}
\begin{equation}   \label{bif17}
\begin{split}
\6z^0_t & = \frac1\eps f^0(\bxm(z^0_t,y^0_t,\eps),z^0_t,y^0_t,\eps) \6t 
+ \frac{\sigma}{\sqrt\eps} F^0(\bxm(z^0_t,y^0_t,\eps),z^0_t,y^0_t,\eps) \6W_t,
\\ 
\6y^0_t & = g(\bxm(z^0_t,y^0_t,\eps),z^0_t,y^0_t,\eps) \6t 
+ \sigma^\prime G(\bxm(z^0_t,y^0_t,\eps),z^0_t,y^0_t,\eps) \6W_t,
\end{split}
\end{equation}
obtained by setting $x^-$ equal to $\bxm(z,y,\eps)$ in~\eqref{bif5}. This
turns out to be true under certain hypotheses on the solutions
of~\eqref{bif17}. Let us fix an initial condition $(z^0_0,y^0_0)\in\cN$,
and call $\z^0_t= (z^0_t,y^0_t)$ the corresponding process. We define the
(random) matrices 
\begin{align}
\label{bif18a}
B(\z^0_t,\eps) &= \left.
\begin{pmatrix}
\sdpar{f^0}z & \sdpar{f^0}y \\
\eps\sdpar gz & \eps\sdpar gy 
\end{pmatrix}
\right\vert_{x=\bxm(z^0_t,y^0_t,\eps),z=z^0_t,y=y^0_t}, \\
\label{bif18b}
C(\z^0_t,\eps) &= \left.
\begin{pmatrix}
\sdpar{f^0}x  \\
\eps\sdpar gx  
\end{pmatrix}
\right\vert_{x=\bxm(z^0_t,y^0_t,\eps),z=z^0_t,y=y^0_t}.
\end{align}
Observe that $C((\hat z,\hat y),0)=0$ because of our choice of coordinates, so
that $\norm{C(\z^0_t,\eps)}$ will be small in a neighbourhood of the origin. 
We denote, for each realization $\z^0(\w)$,  by $\cV_\w$ the principal
solution of 
\begin{equation}
\label{bif19}
\6\z_t(\w)=\frac1\eps B(\z^0_t(\w),\eps)\z_t(\w)\6t. 
\end{equation}
(Note that we may assume that almost all realizations $\z^0(\w)$ are
continuous.) We need to assume the existence of deterministic functions
$\vartheta(t,s)$, $\vartheta_C(t,s)$, and a stopping time
$\tau\leqs\tau_{\cB^-(h)}$ such that 
\begin{equation}
\label{bif20}
\bignorm{\cV_\w(t,s)} \leqs \vartheta(t,s),
\qquad\qquad
\bignorm{\cV_\w(t,s)C(\z^0_s(\w),\eps)} \leqs \vartheta_C(t,s)
\end{equation}
hold for all $s\leqs t\leqs\tau(\w)$ and (almost) all paths
$(\z^0_u(\w))_{u\geqs0}$ of~\eqref{bif17}. Then we define 
\begin{align}
\nonumber
\chi^{(i)}(t) &= \sup_{0\leqs s\leqs t}\frac1\eps \int_0^s 
 \vartheta(s,u)^i \6u,\\
\chi_C^{(i)}(t) &= \sup_{0\leqs s\leqs t}\frac1\eps \int_0^s 
\Bigpar{\sup_{u\leqs v\leqs s} \vartheta_C(s,v)^i} \6u 
\label{pbif21}
\end{align}
for $i=1,2$, and the following result holds. 

\begin{theorem}
\label{thm4}
Assume that there exist constants $\Delta, \vartheta_0>0$ (of order~$1$) such
that $\vartheta(s,u)\leqs\vartheta_0$ and $\vartheta_C(s,u)\leqs\vartheta_0$
whenever $0<s-u\leqs \Delta\eps$. Then there exist constants
$h_0, \kappa_0>0$ such that for all $h\leqs
h_0\brak{\chi^{(1)}(t)\vee\chi_C^{(1)}(t)}^{-1}$ and all initial conditions
$(x^-_0,z^0_0,y^0_0)\in\cB^-(h)$, 
\begin{multline}
\label{bif22}
\qquad\qquad
\Bigprobin{0,(x^-_0,z^0_0,y^0_0)}
{\sup_{0\leqs s\leqs t\wedge\tau} \bignorm{(z_s,y_s)-(z^0_s,y^0_s)} \geqs h}
\\
\leqs \cC_{m,q}(t,\eps) \exp\biggset{-\kappa_0 \frac{h^2}{\sigma^2}
\frac1{\chi_C^{(2)}(t) + h\chi_C^{(1)}(t) + h^2\chi^{(2)}(t)}}, 
\qquad\qquad
\end{multline}
where 
\begin{equation}
\label{bif23}
\cC_{m,q}(t,\eps) = {\it const}\Bigpar{1+\frac t\eps} \e^{(m+q)/4}. 
\end{equation}
\end{theorem}

This result shows that typical solutions of the reduced
system~\eqref{bif17} approximate solutions of the initial
system~\eqref{bif5} to order \smash{$\sigma\chi_C^{(2)}(t)^{1/2} +
\sigma^2\chi_C^{(1)}(t)$}, as long as $\chi^{(1)}(t)\ll
1/\sigma$. Checking the validity of Condition~\eqref{bif20} for a
reasonable stopping time $\tau$ is, of course, not straightforward,
but it depends only on the dynamics of the reduced system, which is
usually easier to analyze. 

\begin{example}   \label{ex_pitchfork}
Assume the reduced equation has the form
\begin{equation}   \label{bif24}
\begin{split}
\6z^0_t & = \frac1\eps \bigbrak{y^0_t z^0_t - (z^0_t)^3} \6t 
+ \frac{\sigma}{\sqrt\eps} \6W_t, \\
\6y^0_t & = 1,
\end{split}
\end{equation}
i.\,e., there is a pitchfork bifurcation at the origin. We fix an
initial time $t_0<0$ and choose an initial condition $(z_0,y_0)$ with
$y_0=t_0$, so that $y^0_t=t$. In~\cite{BG1} we proved that if
$\sigma\leqs\sqrt\eps$, the paths $\set{z_s}_{s\geqs t_0}$ are
concentrated, up to time $\sqrt\eps$, in a strip of width of order
$\sigma/(\abs{y^0}^{1/2}\vee\eps^{1/4})$ around the corresponding
deterministic solution.  

Using for $\tau$ the first-exit time from a set of this form, one finds
that $\chi_C^{(2)}(\sqrt\eps)$ is of order $\sqrt\eps+\sigma^2/\eps$
and that $\chi_C^{(1)}(\sqrt\eps)$ is of order $1+\sigma/\eps^{3/4}$. Thus,
up to time~$\sqrt\eps$, the typical spreading of $z_s$ around reduced
solutions $z^0_s$ is at most of order $\sigma\eps^{1/4} +
\sigma^2/\sqrt\eps$, which is smaller than the spreading of $z^0_s$
around a deterministic solution. Hence the reduced system provides a
good approximation to the full system up to time $\sqrt\eps$. 

For larger times, however, $\chi_C^{(2)}(t)$ grows like 
$\e^{t^2/\eps}$ until the paths leave a neighbourhood of the unstable
equilibrium $z=0$, which typically occurs at a time of order
$\sqrt{\eps\abs{\log\sigma}}$. Thus the spreading is too fast for the reduced
system to provide a good approximation to the dynamics. This shows
that Theorem~\ref{thm4} is not quite sufficient to reduce the problem to a
one-dimensional one, and a more detailed description has to be used for the
region of instability. 
\end{example}


\section{Proofs -- Exit from $\cB(h)$}   \label{sec_exit}

In this section, we consider the SDE
\begin{equation}   \label{exit1}
\begin{split}
\6x_t & = \frac1\eps f(x_t,y_t,\eps) \6t + \frac{\sigma}{\sqrt\eps}
F(x_t,y_t,\eps) \6W_t, \\
\6y_t & = g(x_t,y_t,\eps) \6t + \sigma^\prime G(x_t,y_t,\eps) \6W_t
\end{split}
\end{equation}
under Assumption~\ref{assump_sde}, that is, when starting near a uniformly
asymptotically stable manifold. We denote by $(\xdet_t,\ydet_t)$, with
$\xdet_t=\bx(\ydet_t,\eps)$, the deterministic solution starting in
$\ydet_0=y_0\in \cD_0$. 

The transformation 
\begin{equation}   \label{exit2}
\begin{split}
x_t & = \bx(\ydet_t+\eta_t,\eps) + \xi_t, \\
y_t & = \ydet_t + \eta_t
\end{split}
\end{equation}
yields a system of the form~\eqref{randst4}, which can be written, using
Taylor expansions, as
\begin{align}   \label{exit3}
\6\xi_t & = \frac1\eps \bigbrak{A(\ydet_t,\eps)\xi_t +
  b(\xi_t,\eta_t,t,\eps)} \6t  
+ \frac{\sigma}{\sqrt\eps} \bigbrak{F_0(\ydet_t,\eps) +
  F_1(\xi_t,\eta_t,t,\eps)} \6W_t, \\
\nonumber
\6\eta_t & = \bigbrak{C(\ydet_t,\eps)\xi_t + B(\ydet_t,\eps)\eta_t +
  c(\xi_t,\eta_t,t,\eps)} \6t 
+ \sigma^\prime \bigbrak{G_0(\ydet_t,\eps) + G_1(\xi_t,\eta_t,t,\eps)} \6W_t. 
\end{align}
There are constants $M,M_1$ such that the remainder terms satisfy the bounds
\begin{equation}   \label{exit4}
\begin{split}
\norm{b(\xi,\eta,t,\eps)} 
& \leqs M \bigpar{\norm{\xi}^2 + \norm{\xi}\norm{\eta} +
m\eps\rho^2\sigma^2}, \\ 
\norm{c(\xi,\eta,t,\eps)}
& \leqs M \bigpar{\norm{\xi}^2 + \norm{\eta}^2},\\
\norm{F_1(\xi,\eta,t,\eps)}
& \leqs M_1 \bigpar{\norm{\xi} + \norm{\eta}},\\
\norm{G_1(\xi,\eta,t,\eps)}
& \leqs M_1 \bigpar{\norm{\xi} + \norm{\eta}}
\end{split}
\end{equation}
for all $(\xi,\eta)$ in a compact set and all $t$ such that $\ydet_t\in\cD_0$.
Note that $M$ and $M_1$ may depend on the dimensions~$n$ and~$m$ (see
Remark~\ref{remark_dim}). The term $m\eps\rho^2\sigma^2$ stems from the term
$r(\xi,\eta,t,\eps)$ in~\eqref{randst5}. We shall highlight its $m$-dependence
since it will in general be unavoidable.


\subsection{Timescales of order~$1$}   \label{ssec_timeone}

We first examine the behaviour of $\xi_u$ on an interval $[s,t]$ with
$\Delta = (t-s)/\eps = \orderone{\eps}$. For this purpose, we fix an initial
condition $y_0\in\cD_0$ and assume that $t$ is chosen in such a way that
$\ydet_u\in\cD_0$ for all $u\leqs t$. 

To ease notations, we will not indicate the $\eps$-dependence of $\Xbar(y)$. 
We assume that $\norm{\Xbar(y)} \leqs K_+$ and $\norm{\Xbar(y)^{-1}} \leqs
K_-$ for all $y\in \cD_0$, and define the functions 
\begin{align}   
\nonumber
\Psi(t) & = \frac1\eps \int_0^t \bignorm{\transpose{U(t,u)} \Xbar(\ydet_t)^{-1}
U(t,u)}\6u, \\     
\label{timeone1}
\Phi(t) & = \frac1\eps \int_0^t \Tr\bigbrak{\transpose{U(t,u)}
  \Xbar(\ydet_t)^{-1} U(t,u)}\6u, \\
\nonumber
\Theta(t) & = \frac1\eps \int_0^t \norm{U(t,u)}\6u,
\end{align}
where $U(t,u)$ again denotes the principal solution of $\eps \dot\xi =
A(\ydet_t,\eps)\xi$. Note that the stability of the adiabatic manifold
implies that $\norm{U(t,u)}$ is bounded by a constant times
$\exp\{-K_0(t-u)/\eps\}$, $K_0>0$, for all $t$ and $u\leqs t$. Hence
$\Psi(t)$ and $\Theta(t)$ are of order~$1$, while $\Phi(t)$ is of
order~$n$. In particular, $\Phi(t)\leqs n\Psi(t)$ holds for all times $t$. 

We first concentrate on upper estimates on the probabilities and will deal
with the lower bound in Corollary~\ref{cor_timeone}. Let us remark that
on timescales of order~$1$, we may safely assume that the deviation $\eta_s$
of $y_s$ from its deterministic counterpart remains small. We fix a
deterministic $h_1>0$ and define
\begin{equation}   \label{tau_eta}
\tau_\eta= \inf\bigsetsuch{s>0}{\norm{\eta_s} \geqs h_1}.
\end{equation}
Lemma~\ref{l_eta} below provides an estimate on the tails of the distribution
of $\tau_\eta$. The following proposition estimates the probability that $x_t$
leaves a \lq\lq layer\rq\rq\ similar to $\cB(h)$ during the time interval
$[s,t]$ despite of $\eta_u$ remaining small. Note that in the proposition the
\lq\lq thickness of the layer\rq\rq\ is measured at $\ydet_u$ instead of
$y_u$. 

\begin{prop}   \label{prop_xi}
For all $\alpha\in[0,1)$, all $\gamma\in(0,1/2)$ and all $\mu>0$, 
\begin{align}   
\nonumber
&\sup_{\xi_0 \,\colon\,\pscal{\xi_0}{\Xbar(y_0)^{-1}\xi_0} \leqs \alpha^2h^2
  }\biggprobin{0,(\xi_0,0)}{\sup_{s \leqs u \leqs t \wedge \tau_\eta} 
  \bigpscal{\xi_u}{\Xbar(\ydet_u)^{-1}\xi_u} \geqs h^2} \\
\nonumber
& \ 
{}\leqs \frac{\e^{m\eps\rho^2}}{(1-2\gamma)^{n/2}} 
\exp\biggset{-\gamma\frac{h^2}{\sigma^2}
  \Bigbrak{1-\alpha^2 - M_0\bigpar{\Delta + (1+\mu) h + (h+h_1)\Theta(t)}}} \\
\label{timeone2}
& \  \hphantom{\leqs{}}
{}+\e^{\Phi(t)/4\Psi(t)}
\exp{\biggset{-\frac{h^2}{\sigma^2}
\frac{\mu^2(1-M_0\Delta)}{8M_1^2(\sqrt{K_+}+h_1/h)^2\Psi(t)}}}
\end{align}
holds for all $h<1/\mu$, with a constant $M_0$ depending only on the
linearization $A$ of $f$, $K_+$, $K_-$, $M$, $\norm{F_0}_\infty$, and
on the dimensions~$n$ and~$m$ via $M$. 
\end{prop}

\begin{proof}
The solution of~\eqref{exit3} can be written as 
\begin{align}   \label{timeone3}
\xi_u 
& = U(u)\xi_0 + \frac{\sigma}{\sqrt\eps} 
\int_0^u U(u,v)F_0(\ydet_v,\eps)\6W_v\\
\nonumber
&\hphantom{={}}{}+ \frac{\sigma}{\sqrt\eps} \int_0^u
U(u,v)F_1(\xi_v,\eta_v,v,\eps)\6W_v 
+  \frac{1}{\eps} \int_0^u U(u,v)b(\xi_v,\eta_v,v,\eps)\6v,
\end{align}
where $U(u)=U(u,0)$ as before. Writing $\xi_u = U(u,s)\Upsilon_u$ and defining
\begin{equation}  \label{timeone4}
\tau_\xi = \inf\setsuch{u \geqs 0}{\bigpscal{\xi_u}{\Xbar(\ydet_u)^{-1}\xi_u}
  \geqs h^2}, 
\end{equation}
the probability on the left-hand side of~\eqref{timeone2} can be
rewritten as 
\begin{equation}  \label{timeone5}
P = \biggprobin{0,(\xi_0,0)}{\sup_{s \leqs u\leqs
  t\wedge\tau_\xi\wedge\tau_\eta} \norm{Q(u)\Upsilon_u}\geqs h},
\end{equation}
where $Q(u)=Q_s(u)$ is the symmetric matrix defined by 
\begin{equation}  \label{timeone6}
Q(u)^2 = \transpose{U(u,s)}\Xbar(\ydet_u)^{-1}U(u,s).
\end{equation}
To eliminate the $u$-dependence of $Q$ in~\eqref{timeone5}, we estimate $P$ by 
\begin{equation}  \label{timeone7}
P \leqs \biggprobin{0,(\xi_0,0)}{\sup_{s \leqs u\leqs
  t\wedge\tau_\xi\wedge\tau_\eta} \norm{Q(t)\Upsilon_u}\geqs H},
\end{equation}
where
\begin{equation}  \label{timeone8}
H = h \biggpar{\sup_{s \leqs u\leqs t} \bignorm{Q(u)Q(t)^{-1}}}^{-1}.
\end{equation}
In order to estimate the supremum in~\eqref{timeone8}, we use the fact that
$Q(v)^{-2}$ satisfies the differential equation
\begin{align}  
\nonumber 
\dtot{}{v} Q(v)^{-2} 
& = \frac{1}{\eps} U(s,v) \Bigbrak{-A(\ydet_v)\Xbar(\ydet_v) - \Xbar(\ydet_v)
  \transpose{A(\ydet_v)} + \eps \dtot{}{v} \Xbar(\ydet_v)} \transpose{U(s,v)}
\\
\label{timeone9}
& = \frac{1}{\eps} U(s,v) F_0(\ydet_v,\eps) \transpose{F_0(\ydet_v,\eps)}
\transpose{U(s,v)}, 
\end{align}
and thus
\begin{equation}   \label{timeone10}
Q(u)^{2}Q(t)^{-2} 
= \one + Q(u)^2  \frac1\eps \int_u^t U(s,v)  F_0(\ydet_v,\eps)
\transpose{F_0(\ydet_v,\eps)}\transpose{U(s,v)} \6v 
= \one + \Order{\Delta}.
\end{equation} 
(Recall that $t-u \leqs t-s \leqs \eps\Delta$ in this subsection, which
implies $\norm{U(s,v)}=1+\Order{\Delta}$ and
$\norm{Q(u)^2} \leqs K_-(1+\Order{\Delta})$.) Therefore,
$H=h(1-\Order{\Delta})$. 

We now split $\Upsilon_u$ into three parts, writing 
$\Upsilon_u = \Upsilon^0_u + \Upsilon^1_u + \Upsilon^2_u$, where
\begin{align}  
\nonumber
\Upsilon^0_u 
& = U(s)\xi_0 + \frac{\sigma}{\sqrt\eps} \int_0^u U(s,v)  F_0(\ydet_v,\eps)
\6W_v, \\ 
\label{timeone12}
\Upsilon^1_u 
& =\frac{\sigma}{\sqrt\eps} \int_0^u U(s,v)  F_1(\xi_v,\eta_v,v,\eps) \6W_v, \\
\nonumber 
\Upsilon^2_u 
& =\frac{1}{\eps} \int_0^u U(s,v)  b(\xi_v,\eta_v,v,\eps) \6v,
\end{align}
and estimate $P$ by the sum of the corresponding probabilities
\begin{align}  
\nonumber
P_0 & = \biggprobin{0,(\xi_0,0)}{\sup_{s \leqs u\leqs t}
  \norm{Q(t)\Upsilon^0_u}\geqs H_0}, \\
\label{timeone13}
P_1 & = \biggprobin{0,(\xi_0,0)}{\sup_{s \leqs u\leqs
    t\wedge\tau_\xi\wedge\tau_\eta} \norm{Q(t)\Upsilon^1_u}\geqs H_1}, \\ 
\nonumber
P_2 & =  \biggprobin{0,(\xi_0,0)}{\sup_{s \leqs u\leqs
    t\wedge\tau_\xi\wedge\tau_\eta} \norm{Q(t)\Upsilon^2_u}\geqs H_2}, 
\end{align} 
where $H_0,H_1,H_2$ satisfy $H_0+H_1+H_2=H$. Note that $P_2$ can be
estimated trivially using the fact that
\begin{equation}   \label{timeone15}
\sup_{s \leqs u\leqs t\wedge\tau_\xi\wedge\tau_\eta}\norm{Q(t)\Upsilon^2_u}
\leqs \sqrt{K_-}M(K_+h^2 + \sqrt{K_+}hh_1 +
m\eps\rho^2\sigma^2)(1+\Order{\Delta}) \Theta(t) 
\defby \overbar H_2.
\end{equation}
Now, we choose 
\begin{align}  
\nonumber
H_2 & = 2\overbar H_2, \\
\label{timeone14}
H_1 & = \mu h H, \\
\nonumber
H_0 & = H - H_1 - H_2
\end{align}
for $0<\mu<1/h$, and estimate the remaining probabilities $P_0$ and
$P_1$ by Lemmas~\ref{l_gaussian} and~\ref{l_nongaussian} below. When
estimating $H_0^2$, we may assume $M_0h\Theta(t)<1$, the
bound~\eqref{timeone2} being trivial otherwise. 
\end{proof}

\begin{lemma}   \label{l_gaussian}
Under the hypotheses of Proposition~\ref{prop_xi}, we have for every
$\gamma\in(0,1/2)$, 
\begin{equation}   \label{timeone16}
P_0 
= \biggprobin{0,(\xi_0,0)}{\sup_{s \leqs u\leqs t} 
\norm{Q(t)\Upsilon^0_u}\geqs H_0} 
\leqs \frac1{(1-2\gamma)^{n/2}} \exp\biggset{-\gamma
  \frac{H_0^2-\alpha^2h^2}{\sigma^2}}, 
\end{equation}
holding uniformly for all $\xi_0$ such that
$\pscal{\xi_0}{\Xbar(y_0)^{-1}\xi_0} \leqs \alpha^2h^2$.  
\end{lemma}

\begin{proof}
For every $\widehat\gamma>0$, $(\exp\{\widehat\gamma
\norm{Q(t)\Upsilon^0_u}^2\})_{u \geqs s}$ is a positive
submartingale and, therefore, Doob's submartingale inequality yields
\begin{equation}   \label{timeone18}
P_0 
= \Bigprobin{0,(\xi_0,0)}{\sup_{s \leqs u\leqs t} 
\e^{\widehat\gamma\norm{Q(t)\Upsilon^0_u}^2}
\geqs \e^{\widehat\gamma H_0^2}}
\leqs \e^{-\widehat\gamma H_0^2}
\bigexpecin{0,(\xi_0,0)}{\e^{\widehat\gamma\norm{Q(t)\Upsilon^0_t}^2}}.
\end{equation}
Now, the random variable $Q(t)\Upsilon^0_t$ is Gaussian, with
expectation $E=Q(t)U(s)\xi_0$ and covariance matrix 
\begin{equation}   \label{timeone17}
\Sigma = \frac{\sigma^2}{\eps} Q(t) \biggpar{
\int_0^t U(s,v) F_0(\ydet_v,\eps) \transpose{F_0(\ydet_v,\eps)}
\transpose{U(s,v)} \6v}  \transpose{Q(t)}.
\end{equation}
Thus, using completion of squares to compute the Gaussian integral, we find
\begin{equation}   \label{timeone19}
\bigexpecin{0,(\xi_0,0)}{\e^{\widehat\gamma\norm{Q(t)\Upsilon^0_t}^2}}
=\frac{\e^{\widehat\gamma\pscal{E}{(\one-2\widehat\gamma\Sigma)^{-1}E}}}
{(\det[\one-2\widehat\gamma\Sigma])^{1/2}}.
\end{equation}
By~\eqref{randst13}, we can write
\begin{equation}   \label{timeone20}
\Sigma 
= \sigma^2 Q(t) U(s,t)\bigbrak{\Xbar(\ydet_t) -
U(t)\Xbar(\ydet_0)\transpose{U(t)}}\transpose{U(s,t)}\transpose{Q(t)}
= \sigma^2 \bigbrak{\one - R\transpose{R}},
\end{equation}
where $R=Q(t)U(s)\Xbar(\ydet_0)^{1/2}$, and we have used the fact that
$U(s,t)\Xbar(\ydet_t)\transpose{U(s,t)} = Q(t)^{-2}$. This shows in
particular that
\begin{equation}   \label{timeone21}
\det[\one-2\widehat\gamma\Sigma] \geqs (1-2\widehat\gamma\sigma^2)^n.
\end{equation}
Moreover, since $\norm{R\transpose{R}}=\norm{\transpose{R}R}\in(0,1)$,
we also have
\begin{align} 
\nonumber
\pscal{E}{(\one - 2\widehat\gamma\Sigma)^{-1}E}
& = \pscal{\Xbar(\ydet_0)^{-1/2}\xi_0}{\transpose{R}(\one -
2\widehat\gamma\Sigma)^{-1}R\Xbar(\ydet_0)^{-1/2}\xi_0} \\
\nonumber
&\leqs \alpha^2 h^2 \bignorm{\transpose{R}\bigpar{\one -
2\widehat\gamma\sigma^2\brak{\one-R\transpose{R}}}^{-1}R} \\
\label{timeone22}
&\leqs \alpha^2 h^2 \bigpar{\brak{1 - 2\widehat\gamma\sigma^2}
\norm{\transpose{R}R}^{-1} + 2\widehat\gamma\sigma^2}^{-1}
\leqs \alpha^2 h^2
\end{align}
for all $\xi_0$ satisfying $\pscal{\xi_0}{\Xbar(y_0)^{-1}\xi_0} \leqs
\alpha^2h^2$. Now, \eqref{timeone16} follows from~\eqref{timeone19}
by choosing $\widehat\gamma=\gamma/\sigma^2$. 
\end{proof}

\begin{lemma}   \label{l_nongaussian}
Under the hypotheses of Proposition~\ref{prop_xi},
\begin{equation}   \label{timeone23}
P_1 = \biggprobin{0,(\xi_0,0)}{\sup_{s \leqs u\leqs
  t\wedge\tau_\xi\wedge\tau_\eta} \norm{Q(t)\Upsilon^1_u}\geqs H_1}
  \!\leqs\! 
\exp{\biggset{-\frac{\bigpar{H_1^2-\sigma^2 M_1^2(\sqrt{K_+}h+h_1)^2\Phi(t)}^2}
  {8\sigma^2M_1^2(\sqrt{K_+}h +h_1)^2 H_1^2\Psi(t)}}}
\end{equation}
holds uniformly for all $\xi_0$ such that
$\pscal{\xi_0}{\Xbar(y_0)^{-1}\xi_0} \leqs h^2$.
\end{lemma}

\begin{proof}
Let $\tau$ denote the stopping time
\begin{equation}   \label{timeone23.5}
\tau=\tau_\xi \wedge \tau_\eta \wedge 
  \inf\setsuch{u\geqs 0}{\norm{Q(t)\Upsilon^1_u} \geqs H_1},
\end{equation}
and define, for a given $\gamma_1$, the stochastic process 
\begin{equation}   \label{timeone24}
\Xi_u = \e^{\gamma_1\norm{Q(t)\Upsilon^1_u}^2}.
\end{equation}
$(\Xi_u)_u$ being a positive submartingale, another application of
Doob's submartingale inequality yields 
\begin{equation}   \label{timeone25}
P_1 \leqs \e^{-\gamma_1H_1^2}
\bigexpecin{0,(\xi_0,0)}{\Xi_{t \wedge \tau}}. 
\end{equation}
It\^o's formula (together with the fact that
$\transpose{(\6W_u)}\transpose{R}R\6W_u = \Tr(\transpose{R}R)\6u$ for
any matrix $R\in\R^{n\times k}$) shows that $\Xi_u$ obeys the SDE
\begin{equation}   \label{timeone26}
\6\Xi_u 
= 2\gamma_1 \frac{\sigma}{\sqrt\eps} \Xi_u
\transpose{(\Upsilon^1_u)} Q(t)^2 U(s,u) F_1(\xi_u,\eta_u,u,\eps)
\6W_u 
+ \gamma_1 \frac{\sigma^2}{\eps} \Xi_u \Tr\bigbrak{\transpose{R_1}R_1
+ 2\gamma_1\transpose{R_2}R_2}\6u,
\end{equation} 
where
\begin{equation}   \label{timeone27}
\begin{split}
R_1 & = Q(t) U(s,u)F_1(\xi_u,\eta_u,u,\eps), \\
R_2 & = \transpose{(\Upsilon^1_u)} Q(t)^2 U(s,u) F_1(\xi_u,\eta_u,u,\eps).
\end{split}
\end{equation} 
The first term in the trace can be estimated as
\begin{align}  
\nonumber
\Tr\bigbrak{\transpose{R_1}R_1}
= \Tr\bigbrak{R_1\transpose{R_1}}
& \leqs M_1^2 \bigpar{\norm{\xi_u} + \norm{\eta_u}}^2
\Tr\bigbrak{\transpose{Q(t)}U(s,u)\transpose{U(s,u)}Q(t)} \\ 
\label{timeone28}
& \leqs M_1^2 \bigpar{\norm{\xi_u} + \norm{\eta_u}}^2
\Tr\bigbrak{\transpose{U(t,u)}\Xbar(\ydet_t)^{-1}U(t,u)},
\end{align}
while the second term satisfies the bound
\begin{align}  
\nonumber
\Tr\bigbrak{\transpose{R_2}R_2}
& = \bignorm{\transpose{F_1(\xi_u,\eta_u,u,\eps)} \transpose{U(s,u)}
Q(t)^2 \Upsilon^1_u}^2 \\
\nonumber
& \leqs M_1^2 \bigpar{\norm{\xi_u} + \norm{\eta_u}}^2
 \bignorm{\transpose{U(s,u)} Q(t)}^2 \norm{Q(t)\Upsilon^1_u}^2\\
\label{timeone29}
& = M_1^2 \bigpar{\norm{\xi_u} + \norm{\eta_u}}^2
\bignorm{\transpose{U(t,u)}\Xbar(\ydet_t)^{-1}
U(t,u)} \norm{Q(t)\Upsilon^1_u}^2.
\end{align}
Using the fact that $\norm{\xi_u} \leqs \sqrt{K_+} h$,
$\norm{\eta_u}\leqs h_1$ and $\norm{Q(t)\Upsilon^1_u} \leqs H_1$
hold for all $0 \leqs u \leqs t\wedge\tau$, we obtain  
\begin{align}   
\nonumber
&\bigexpecin{0,(\xi_0,0)}{\Xi_{u\wedge\tau}}
\\ 
\nonumber
& \qquad \leqs 1 + \gamma_1\frac{\sigma^2}{\eps} M_1^2
\bigpar{\sqrt{K_+}h+h_1}^2 
\int_0^u \bigexpecin{0,(\xi_0,0)}{\Xi_{v\wedge\tau}} \\ 
\label{timeone30}
&\qquad \phantom{\leqs{}}
{}\times\Bigbrak{\Tr\bigbrak{\transpose{U(t,v)}\Xbar(\ydet_t)^{-1}U(t,v)}
+  2\gamma_1 H_1^2 \bignorm{\transpose{U(t,v)}\Xbar(\ydet_t)^{-1}
U(t,v)}}\6v,
\end{align}
and Gronwall's inequality yields
\begin{equation}    \label{timeone31a}
\bigexpecin{0,(\xi_0,0)}{\Xi_{t\wedge\tau}}
\leqs \exp\Bigset{\gamma_1\sigma^2M_1^2 \bigpar{\sqrt{K_+}h+h_1}^2
\bigbrak{\Phi(t) + 2\gamma_1 H_1^2 \Psi(t)}}.
\end{equation}
Now, \eqref{timeone25} implies 
\begin{equation}    \label{timeone31b}
P_1 \leqs \exp\Bigset{-\gamma_1\bigpar{H_1^2 - \sigma^2
    M_1^2(\sqrt{K_+}h+h_1)^2 
\Phi(t)} + 2\gamma_1^2 \sigma^2 M_1^2(\sqrt{K_+}h+h_1)^2 H_1^2\Psi(t)}, 
\end{equation}
and~\eqref{timeone23} follows by optimizing over $\gamma_1$.
\end{proof}

Proposition~\ref{prop_xi} allows to control the first-exit time of $(x_t,y_t)$ 
from $\cB(h)$, provided $\eta_s = y_s - \ydet_s$ remains small. In order to
complete the proof of Part~(a) of Theorem~\ref{thm1} we need to control the
tails of the distribution of $\tau_\eta$. The following lemma provides a rough
{\it a priori\/} estimate which is sufficient for the time being. We will
provide more precise estimates in the next section. 

Recall the notations $V(u,v)$ for the principal solution of $\dot\eta =
B(\ydet_u,\eps)\eta$, and 
\begin{align}    \label{timeone43a}
\chi^{(1)}(t) & = \sup_{0\leqs s \leqs t} \int_0^s \Bigpar{\sup_{u\leqs v\leqs
    s} \norm{V(s,v)}} \6u,  \\
\label{timeone43b}
\chi^{(2)}(t) & = \sup_{0\leqs s \leqs t} \int_0^s \Bigpar{\sup_{u\leqs v\leqs
    s} \norm{V(s,v)}^2} \6u.
\end{align}
from Subsection~\ref{ssec_randst}.

\begin{lemma}   \label{l_eta}
There exists a constant $c_\eta>0$ such that for all choices of $t>0$ and
$h_1>0$ satisfying $\ydet_s\in\cD_0$ for all $s\leqs t$ and $h_1\leqs
c_\eta\chi^{(1)}(t)^{-1}$,
\begin{align}
\nonumber
&\sup_{\xi_0 \,\colon\,\pscal{\xi_0}{\Xbar(y_0)^{-1}\xi_0} \leqs h^2}
\biggprobin{0,(\xi_0,0)}{\sup_{0 \leqs u \leqs t \wedge \tau_{\cB(h)}}  
\norm{\eta_u} \geqs h_1} \\
\nonumber
&\qquad\qquad \leqs 2\biggintpartplus{\frac{t}{\Delta\eps}}\e^{m/4}
\exp\biggset{\frac{m\eps\rho^2}{(\rho^2+\eps)\chi^{(2)}(t)}} \\
\label{timeone42}
&\qquad\qquad \phantom{{}\leqs{}}
{}\times  \exp\biggset{-\kappa_0
\frac{h_1^2(1-\Order{\Delta\eps})}{\sigma^2(\rho^2+\eps)\chi^{(2)}(t)}
\biggbrak{1-M_0^\prime\,\chi^{(1)}(t)\,h_1\Bigpar{1+K_+\frac{h^2}{h_1^2}}}}, 
\end{align}
where $\kappa_0>0$ is a constant depending only on $\norm{\widehat
F}_\infty$, $\norm{\widehat G}_\infty$, $\norm{C}_\infty$ and $U$, while
the constant $M_0^\prime$ depends only on $M$, $\norm{C}_\infty$ and $U$. 
Note that $c_\eta$ may depend on the dimensions~$n$ and~$m$ via $M$.
\end{lemma}

In the sequel, we will typically choose $h_1\gg\sigma$, so that the
prefactor becomes negligible.

\begin{proof}[\sc Proof of Lemma~\ref{l_eta}]
We first consider a time interval $[s,t]$ with $t-s=\Delta\eps$. Let
$u\in[s,t]$ and recall the defining SDE~\eqref{exit3} for $\eta_u$. Its
solution can be split into four parts, $\eta_u = \eta^0_u + \eta^1_u +
\eta^2_u + \eta^3_u$, where  
\begin{equation}
\begin{split}
\eta^0_u &= \sigma^\prime \int_0^u V(u,v) \widehat
G(\xi_v,\eta_v,v,\eps)\6W_v, \\
\eta^1_u &= \frac\sigma{\sqrt\eps} \int_0^u S(u,v) \widehat
F(\xi_v,\eta_v,v,\eps)\6W_v, \\
\eta^2_u &= \int_0^u V(u,v) c(\xi_v,\eta_v,v,\eps)\6v, \\
\label{timeone44}
\eta^3_u &= \frac1\eps\int_0^u S(u,v) b(\xi_v,\eta_v,v,\eps)\6v,
\end{split}
\end{equation}
with
\begin{equation}    \label{timeone45}
S(u,v) = \int_v^u V(u,w) C(\ydet_w,\eps) U(w,v) \6w.
\end{equation}
Let $\tau=\tau_{\cB(h)}\wedge\tau_\eta$. It follows immediately from the
definitions of $\tau_{\cB(h)}$, $\tau_\eta$ and the bounds~\eqref{exit4} that 
\begin{equation}
\begin{split}
\bignorm{\eta^2_{u\wedge\tau}} &\leqs M (1+\Order{\Delta\eps}) \chi^{(1)}(t)
(K_+h^2 + h_1^2), \\ 
\label{timeone46}
\bignorm{\eta^3_{u\wedge\tau}} 
&\leqs M^\prime \chi^{(1)}(t) (K_+h^2 + \sqrt{K_+}hh_1 + m\eps\rho^2\sigma^2) 
\end{split}
\end{equation}
for all $u\in[s,t]$. Here $M^\prime$ depends only on $M$, $U$ and
$\norm{C}_\infty$. Furthermore, using similar ideas as in the proof of
Lemma~\ref{l_nongaussian}, it is straightforward to establish for all $H_0,
H_1> 0$ that 
\begin{equation}
\begin{split}\biggprobin{0,(\xi_0,0)}{\sup_{s \leqs u\leqs t\wedge\tau}
\norm{\eta^0_u}\geqs H_0}
&\leqs \e^{m/4} \exp\biggset{-\frac{H_0^2(1-\Order{\Delta\eps})}
{8(\sigma^\prime)^2\norm{\widehat G}_\infty^2\chi^{(2)}(t)}}, \\
\label{timeone47}
\biggprobin{0,(\xi_0,0)}{\sup_{s \leqs u\leqs t\wedge\tau}
\norm{\eta^1_u}\geqs H_1}
&\leqs \e^{m/4} \exp\biggset{-\frac{H_1^2(1-\Order{\Delta\eps})}
{8\sigma^2\eps c_S\norm{\widehat F}_\infty^2\chi^{(2)}(t)}}, 
\end{split}
\end{equation}
where $c_S$ is a constant depending only on $S$. Then the local analogue of 
estimate~\eqref{timeone42} (without the $t$-dependent prefactor) is obtained
by taking, for instance, 
$H_0=H_1= \frac12 h_1 - 2(M+M^\prime)\chi^{(1)}(t)(K_+h^2 + h_1^2 +
m\eps\rho^2\sigma^2)$, and using $h_1\leqs c_\eta\chi^{(1)}(t)^{-1}$,
where we may choose $c_\eta\leqs 1/(2M_0^\prime)$. 

It remains to extend~\eqref{timeone42} to a general time interval
$[0,t]$ for $t$ of order~$1$. For this purpose, we choose a partition 
$0=u_0<u_1<\dots<u_K=t$ of $[0,t]$, satisfying $u_k=k\Delta\eps$ for $0\leqs
k<K=\intpartplus{t/(\Delta\eps)}$. Applying the local version
of~\eqref{timeone42} to each interval $[u_k,u_{k+1}]$ and using the
monotonicity of $\chi^{(2)}(u)$, the claimed estimate follows from
\begin{equation}   \label{timeone47.5}
\vphantom{\sum_i} \Bigprobin{0,(\xi_0,0)}{\sup_{0 \leqs u \leqs t \wedge
    \tau_{\cB(h)}} \norm{\eta_u} \geqs h_1} 
\leqs \sum_{k=0}^{K-1}
\Bigprobin{0,(\xi_0,0)}{\sup_{u_k \leqs u \leqs u_{k+1} \wedge
    \tau_{\cB(h)}} \norm{\eta_u} \geqs h_1}.
\end{equation}
\end{proof}

We will now show that Proposition~\ref{prop_xi} and Lemma~\ref{l_eta} together
are sufficient to prove Parts~(a) and~(b) of Theorem~\ref{thm1} on a timescale
of order~$1$. We continue to assume that $y_0\in\cD_0$ but we will no longer
assume that $\ydet_u\in\cD_0$ automatically holds for all $u\leqs t$. Instead,
we will employ Lemma~\ref{l_eta} to compare $y_u\in\cD_0$ and $\ydet_u$,
taking advantage of the fact that on timescales of order~$1$, $\eta_t$ is
likely to remain small. Note that if the uniform-hyperbolicity
Assumption~\ref{assump_sde} holds for $\cD_0$, then there exists a $\delta>0$
of order~$1$ such that the $\delta$-neighbourhood $\cD_0^+(\delta)$ 
also satisfies this assumption. We introduce the first-exit time $\taudet$ of
the deterministic process $\ydet_u$ from $\cD_0^+(\delta)$ as
\begin{equation}   \label{timeone47.6}
\taudet = \inf\setsuch{u\geqs0}{\ydet_u\not\in\cD_0^+(\delta)}
\end{equation}
and remark in passing that $\tau_{\cB(h)} \wedge \tau_{\eta} \leqs \taudet$
holds whenever $h_1\leqs\delta$.

\begin{cor}    \label{cor_timeone}
Fix a time $t>0$ and $h>0$ in such a way that 
$h\leqs c_1 \chi^{(1)}(t\wedge\taudet)^{-1}$ for a sufficiently small constant
$c_1 > 0$ and $\chi^{(2)}(t\wedge\taudet) \leqs (\rho^2+\eps)^{-1}$.  
Then for any $\alpha\in[0,1)$, any $\gamma\in(0,1/2)$ and any sufficiently
small $\Delta$,  
\begin{equation}   \label{timeone48}
 C^-_{n,m}(t,\eps) \e^{-\kappa^-(0) h^2/\sigma^2}
\leqs \bigprobin{0,(\xi_0,0)}{\tau_{\cB(h)} < t} 
\leqs C^+_{n,m,\gamma}(t,\eps) \e^{-\kappa^+(\alpha) h^2/\sigma^2}
\end{equation}
holds uniformly for all $\xi_0$ satisfying $\pscal{\xi_0}{\Xbar(y_0)^{-1}\xi_0}
\leqs \alpha^2h^2$. Here
\begin{align}   
\label{timeone49a}
\vphantom{\biggpar{}}
\kappa^+(\alpha) &= \gamma \bigbrak{1 - \alpha^2 - \Order{\Delta}
-\Order{m\eps\rho^2} - \Order{(1+\chi^{(1)}(t\wedge\taudet))h}}, \\
\label{timeone49b}
\kappa^-(0) &= \frac12\bigbrak{1+\Order{h}+\bigOrder{\e^{-K_0t/\eps}}}, \\
\label{timeone49c}
C^+_{n,m,\gamma}(t,\eps) &= 
\biggintpartplus{\frac{t}{\Delta\eps}} 
\biggbrak{\frac1{(1-2\gamma)^{n/2}}
+ \bigpar{\e^{n/4} + 2\e^{m/4}}\e^{-\kappa^+(0) h^2/\sigma^2}}, \\
\label{timeone49d}
C^-_{n,m}(t,\eps) &= 
\biggpar{\sqrt{\frac2\pi}\frac h\sigma \wedge
1}\e^{-\Order{m\eps\rho^2}} \\
\nonumber
&\phantom{{}={}} {} - \biggpar{\e^{n/4} +
4\biggintpartplus{\frac{t}{\Delta\eps}} \e^{m/4}}  
\e^{-\frac{h^2}{2\sigma^2}\brak{1-\Order{\e^{-K_0t/\eps}}
-\Order{m\eps\rho^2} - \Order{(1+\chi^{(1)}(t\wedge\taudet))h}}}.
\end{align}
\end{cor}

\begin{proof}
We first establish the upper bound. Fix an initial condition $(\xi_0,0)$
satisfying 
$\pscal{\xi_0}{\Xbar(y_0)^{-1}\xi_0} \leqs \alpha^2h^2$, and observe that 
\begin{align}   \label{timeone52}
\bigprobin{0,(\xi_0,0)}{\tau_{\cB(h)} < t} 
& \leqs  \bigprobin{0,(\xi_0,0)}{\tau_{\cB(h)} < t\wedge\tau_\eta} 
+ \bigprobin{0,(\xi_0,0)}{\tau_\eta < t\wedge\tau_{\cB(h)}} \\
\nonumber
& =  
\bigprobin{0,(\xi_0,0)}{\tau_{\cB(h)} < t\wedge\taudet\wedge\tau_\eta} 
+ \bigprobin{0,(\xi_0,0)}{\tau_\eta < t\wedge\taudet\wedge\tau_{\cB(h)}}.
\end{align}
To estimate the first term on the right-hand side, we again
introduce a partition $0=u_0<u_1<\dots<u_K=t$ of the time interval $[0,t]$,
defined by $u_k=k\Delta\eps$ for $0\leqs
k<K=\intpartplus{t/(\Delta\eps)}$. Thus we obtain
\begin{equation}   \label{timeone52.5}
\bigprobin{0,(\xi_0,0)}{\tau_{\cB(h)} < t\wedge\taudet\wedge\tau_\eta}
\leqs  \sum_{k=1}^K \bigprobin{0,(\xi_0,0)}{u_{k-1} \leqs 
   \tau_{\cB(h)} <  u_k \wedge\taudet\wedge \tau_\eta}.
\end{equation}
Before we estimate the summands on the right-hand side of~\eqref{timeone52.5}, 
note that by the boundedness assumption on $\norm{\Xbar(y)}$ and
$\norm{\Xbar^{-1}(y)}$, we have $\Xbar(y_u)^{-1} = \Xbar(\ydet_u)^{-1} +
\Order{h_1}$ for $u\leqs\taudet\wedge\tau_\eta$. Thus the bound obtained in
Proposition~\ref{prop_xi} can also be applied to estimate first-exit times
from $\cB(h)$ itself: 
\begin{align}
\nonumber
&\bigprobin{0,(\xi_0,0)}{u_{k-1} \leqs\tau_{\cB(h)} < u_k
  \wedge\taudet\wedge\tau_\eta} \\  
\label{timeone50}
&\qquad{}\leqs \Bigprobin{0,(\xi_0,0)}
{\sup_{u_{k-1}\leqs u<u_k \wedge\taudet\wedge\tau_\eta}
  \bigpscal{\xi_u}{\Xbar(\ydet_u)^{-1}\xi_u} 
\geqs h^2(1-\Order{h_1})},
\end{align}
while the second term on the right-hand side of~\eqref{timeone52} can
be estimated directly by Lemma~\ref{l_eta}. Choosing 
\begin{equation}   \label{timeone50.5}
\mu^2 =
8M_1^2\bigbrak{\sqrt{K_+}+h_1/(h(1-\Order{h_1}))}^2\Psi(t\wedge\taudet)
/\bigbrak{1-\Order{h_1}-M_0\Delta} 
\end{equation}
and $h_1=h/\sqrt{\kappa_0}$ in the resulting expression, we see that
the Gaussian part of $\xi_t$ gives the major contribution to the
probability. Thus we obtain that the probability in~\eqref{timeone52}
is bounded by  
\begin{align}
\nonumber
&\biggintpartplus{\frac{t}{\Delta\eps}} \biggbrak{
\frac{\e^{m\eps\rho^2}}{(1-2\gamma)^{n/2}} 
\exp\Bigset{-\gamma\frac{h^2}{\sigma^2}
\bigbrak{1-\alpha^2-\Order{\Delta}-\Order{h}}} 
+ \e^{n/4}\e^{-h^2/\sigma^2} \\
\label{timeone51}
&\qquad\quad {}+ 2\e^{m/4}
\exp\biggset{-\frac{h^2
    (1-\Order{\chi^{(1)}(t\wedge\taudet)h}-\Order{\Delta\eps} -
\Order{m\eps\rho^2})}{\sigma^2 (\rho^2+\eps) \chi^{(2)}(t\wedge\taudet)}}}, 
\end{align}
where we have used the fact that $\Phi(t)\leqs n\Psi(t)$, while $\Psi(t)$ and
$\Theta(t)$ are at most of order~1. The prefactor $\e^{m\eps\rho^2}$ can be
absorbed into the error term $\Order{m\eps\rho^2}$ in the exponent. This
completes the proof of the upper bound in~\eqref{timeone48}.

The lower bound is a consequence of the fact that the Gaussian part of $\xi_t$
gives the major contribution to the probability in~\eqref{timeone48}. To check
this, we split the probability as follows:
\begin{align}
\nonumber
& \bigprobin{0,(\xi_0,0)}{\tau_{\cB(h)} < t} \\
\nonumber
& \qquad{}\geqs \bigprobin{0,(\xi_0,0)}{\tilde\tau_\xi < t, \tau_\eta \geqs t} 
+ \bigprobin{0,(\xi_0,0)}{\tau_{\cB(h)} < t, \tau_\eta < t} \\
\nonumber
& \qquad{}=  \bigprobin{0,(\xi_0,0)}{\tilde\tau_\xi < t \wedge \tau_\eta} 
- \bigprobin{0,(\xi_0,0)}{\tilde\tau_\xi < \tau_\eta < t}
 + \bigprobin{0,(\xi_0,0)}{\tau_{\cB(h)} < t, \tau_\eta < t} \\
\label{timeone53} 
& \qquad{}\geqs \bigprobin{0,(\xi_0,0)}{\tilde\tau_\xi < t \wedge \tau_\eta} 
- \bigprobin{0,(\xi_0,0)}{\tau_\eta < t \wedge \tau_{\cB(h)}},
\end{align}
where
\begin{equation}   \label{timeone54}
\tilde\tau_\xi = \inf\setsuch{u \geqs
  0}{\bigpscal{\xi_u}{\Xbar(\ydet_u)^{-1}\xi_u} \geqs h^2(1+\Order{h_1})},
\end{equation}
and the $\Order{h_1}$-term stems from estimating $\Xbar(y_u)^{-1}$ by
$\Xbar(\ydet_u)^{-1}$ as in~\eqref{timeone50}. 
The first term on the last line of~\eqref{timeone53} can be estimated
as in the proof of Proposition~\ref{prop_xi}: A lower bound is
obtained trivially by considering the endpoint instead of the whole
path, and instead of applying Lemma~\ref{l_gaussian}, the Gaussian
contribution can be estimated below by a straightforward
calculation. The non-Gaussian parts are estimated {\it above\/} as
before and are of smaller order. Finally, we need an upper 
bound for the probability that $\tau_\eta < t \wedge \tau_{\cB(h)}$,
which can be obtained from Lemma~\ref{l_eta}.  
\end{proof}


\subsection{Longer timescales}   \label{ssec_timetwo}

Corollary~\ref{cor_timeone} describes the dynamics on a timescale of order
$1$, or even on a slightly longer timescale if $\chi^{(1)}(t)$,
$\chi^{(2)}(t)$ do not grow too fast. It may happen, however, that $\ydet_t$
remains in $\cD_0$ for all positive times (e.\,g.\ when $\cD_0$ is positively
invariant under the reduced deterministic flow). In such a case, one would
expect the vast majority of paths to remain concentrated in $\cB(h)$ for a
rather long period of time.   

The approach used in Subsection~\ref{ssec_timeone} fails to control the
dynamics on timescales on which $\chi^{(i)}(t)\gg1$, because it uses in an
essential way the fact that $\eta_t = y_t-\ydet_t$ remains small. Our
strategy in order to describe the paths on longer timescales is to compare
them to different deterministic solutions on time intervals $[0,T]$,
$[T,2T]$, \dots, where $T$ is a possibly large constant such that
Corollary~\ref{cor_timeone} holds on time intervals of length $T$, provided
$y_t$ remains in $\cD_0$. Essential ingredients for this approach are the
Markov property and the following technical lemma, which is based on
integration by parts.  

\begin{lemma}   \label{l_intbyparts}
Fix constants $s_1\leqs s_2$ in $[0,\infty]$, and assume we are given two
continuously differentiable functions 
\begin{itemiz}
\item   $\varphi: [0,\infty) \to [0,\infty)$, which is monotonously
increasing and satisfies $\varphi(s_2)=1$,
\item   $\varphi_0: [0,\infty) \to \R$ which satisfies $\varphi_0(s)\leqs 0$
  for all $s\leqs s_1$. 
\end{itemiz}
Let $X\geqs0$ be a random variable such that 
$\prob{X<s} \geqs \varphi_0(s)$ for all $s\geqs0$. 
Then we have, for all $t\geqs0$, 
\begin{equation}   \label{timetwo1}
\bigexpec{\indicator{[0,t)}(X)\hat\varphi(X)} 
\leqs \hat\varphi(t)\prob{X<t} - \int_{s_1\wedge t}^{s_2\wedge t}
\varphi^\prime(s)\varphi_0(s)\6s,
\end{equation}
where $\hat\varphi(s)=\varphi(s)\wedge1$. 
\end{lemma}

We omit the proof of this result, which is rather standard. See, for
instance, \cite[Lemma~A.1]{BG1} for a very similar result. 

When applying the preceding lemma, we will also need an estimate on the
probability that $\pscal{\xi_T}{\Xbar(y_T)^{-1}\xi_T}$ exceeds $h^2$.
Corollary~\ref{cor_timeone} provides, of course, such an estimate, but
since it applies to the whole path, it does not give optimal bounds for the
endpoint. An improved bound is given by the following lemma. Recall the
definition of the first-exit time $\tau_{\cD_0}$ of $y_t$ from $\cD_0$
from~\eqref{randst17}. 

\begin{lemma}   \label{l_endpoint}
If $T$ and $h$ satisfy $h\leqs c_1 \chi^{(1)}(T\wedge\taudet)^{-1}$ and
$\chi^{(2)}(T\wedge\taudet) \leqs (\rho^2+\eps)^{-1}$, we have, for
every $\gamma\in(0,1/2)$, 
\begin{equation}   \label{timetwo2}
\sup_{\xi_0 \,\colon\,\pscal{\xi_0}{\Xbar(y_0)^{-1}\xi_0} \leqs h^2} 
\Bigprobin{0,(\xi_0,0)}{\pscal{\xi_T}{\Xbar(y_T)^{-1}\xi_T}\geqs h^2,
  \tau_{\cD_0}\geqs T}  
\leqs \widehat C_{n,m,\gamma}(T,\eps) \e^{-\kappa^\prime h^2/\sigma^2},
\end{equation}
where 
\begin{align}   
\label{timetwo3a}
\kappa^\prime &= \gamma \bigbrak{1 - \Order{\Delta} - \Order{h} -
\bigOrder{\e^{-2K_0T/\eps}/(1-2\gamma)}},
\\
\label{timetwo3b}
\widehat C_{n,m,\gamma}(T,\eps) &= 
\frac{\e^{m\eps\rho^2}}{(1-2\gamma)^{n/2}} 
+ 4 C^+_{n,m,\gamma}(T,\eps) \e^{-2\kappa^+(0) h^2/\sigma^2}.
\end{align}
\end{lemma}
\begin{proof}
We decompose $\xi_t$ as $\xi_t=\xi^0_t + \xi^1_t + \xi^2_t$, where 
\begin{align}
\nonumber
\xi^0_t &= U(t)\xi_0 + \frac\sigma{\sqrt\eps} \int_0^t U(t,u) 
F_0(\ydet_u,\eps)\6W_u, \\
\label{timetwo4}
\xi^1_t &= \frac\sigma{\sqrt\eps} \int_0^t U(t,u) 
F_1(\xi_u,\eta_u,u,\eps)\6W_u, \\
\nonumber
\xi^2_t &= \frac1\eps \int_0^t U(t,u) b(\xi_u,\eta_u,u,\eps)\6u,
\end{align}
and introduce the notations $\tilde\tau_\xi$ and $\tilde\tau_\eta$ for the
stopping times which are defined like $\tau_\xi$ and $\tau_\eta$
in~\eqref{timeone4} and~\eqref{tau_eta}, but with $h$ and $h_1$ replaced by
$2h$ and $2h_1$, respectively. The probability in~\eqref{timetwo2} is bounded
by 
\begin{equation}   \label{timetwo5}
\Bigprobin{0,(\xi_0,0)}{\pscal{\xi_T}{\Xbar(\ydet_T)^{-1}\xi_T}\geqs
h^2(1-\Order{h_1}), \tilde\tau_\eta>T} 
+ \Bigprobin{0,(\xi_0,0)}{\tilde\tau_\eta\leqs T}.
\end{equation}
Let $H^2 = h^2(1-\Order{h_1})$. As in the proof of Proposition~\ref{prop_xi},
the first term can be further decomposed as 
\begin{align}
\nonumber
&\Bigprobin{0,(\xi_0,0)}{\pscal{\xi_T}{\Xbar(\ydet_T)^{-1}\xi_T}\geqs H^2,
    \tilde\tau_\eta>T} \\
\nonumber
&\qquad \leqs 
\Bigprobin{0,(\xi_0,0)}{\bignorm{\Xbar(\ydet_T)^{-1/2}\xi^0_T} \geqs H_0} 
+ \Bigprobin{0,(\xi_0,0)}{ \tilde\tau_\eta >T, \tilde\tau_\xi \leqs T}  \\
\nonumber
&\qquad\phantom{\leqs{}}
{}+ \Bigprobin{0,(\xi_0,0)}{\bignorm{\Xbar(\ydet_T)^{-1/2}\xi^1_T} \geqs H_1,
  \tilde\tau_\eta >T, \tilde\tau_\xi >T} \\ 
\label{timetwo6}
&\qquad\phantom{\leqs{}}
{}+ \Bigprobin{0,(\xi_0,0)}{\bignorm{\Xbar(\ydet_T)^{-1/2}\xi^2_T} \geqs H_2,
 \tilde\tau_\eta >T, \tilde\tau_\xi >T},
\end{align}
where we choose $H_1$, $H_2$ twice as large as in the proof of
Proposition~\ref{prop_xi}, while $H_0=H-H_1-H_2$. 

The first term on the right-hand side can be estimated as
in~Lemma~\ref{l_gaussian}, with the difference that, the expectation
of $\xi^0_T$ being exponentially small in $T/\eps$, it leads only to a
correction of order $\e^{-2K_0T/\eps}/(1-2\gamma)$ in the exponent. The
second and the third term can be estimated by Corollary~\ref{cor_timeone} and
Lemma~\ref{l_nongaussian}, the only difference lying in a larger absolute
value of the exponent, because we enlarged $h$ and $h_1$. The last term
vanishes by our choice of $H_2$. Finally, the second term
in~\eqref{timetwo5} can be estimated by splitting according to the
value of $\tau_{\cB(2h)}$ and applying Lemma~\ref{l_eta} and
Corollary~\ref{cor_timeone}.  
\end{proof}

We are now ready to establish an improved estimate on the distribution of
$\tau_{\cB(h)}$. As we will restart the process $\ydet_t$ whenever $t$ is a
multiple of $T$, we need the assumptions made in the previous section to hold
uniformly in the initial condition $y_0\in\cD_0$. Therefore we will introduce
replacements for some of the notations introduced before. Note that
$\chi^{(1)}(t)=\chi_{y_0}^{(1)}(t)$ and $\chi^{(2)}(t)=\chi_{y_0}^{(2)}(t)$
depend on $y_0$ via the principal solution $V$. Also $\taudet=\taudet(y_0)$
naturally depends on $y_0$. We define 
\begin{align}
\label{chihat1}
\widehat\chi^{(1)}(t) & = \sup_{y_0\in\cD_0}
\chi_{y_0}^{(1)}\bigpar{t\wedge\taudet(y_0)},  \\ 
\label{chihat2}
\widehat\chi^{(2)}(t) & = \sup_{y_0\in\cD_0}
\chi_{y_0}^{(2)}\bigpar{t\wedge\taudet(y_0)}.
\end{align}
In the same spirit, the $\chi^{(i)}(T)$-dependent $\Order{\cdot}$-terms in the
definitions of $\kappa^+(\alpha)$, $\kappa'$ and the
prefactors like $C^+_{n,m,\gamma}(T,\eps)$ are modified.

We fix a time $T$ of order~$1$ satisfying $\widehat\chi^{(2)}(T) \leqs
(\rho^2+\eps)^{-1}$. $T$ is chosen in such a way that whenever $h\leqs c_1
\widehat\chi^{(1)}(T)^{-1}$, Corollary~\ref{cor_timeone} (and
Lemma~\ref{l_endpoint}) apply. Note that larger $T$ would be possible unless
$\rho$ is of order~$1$, but for larger $T$ the constraint on $h$ becomes more
restrictive which is not desirable. Having chosen $T$, we define the
probabilities 
\begin{align}
\label{timetwo8a}
P_k(h) &= \bigprobin{0,(0,0)}{\tau_{\cB(h)}<kT \wedge \tau_{\cD_0}}, \\
\label{timetwo8b}
Q_k(h) &= \bigprobin{0,(0,0)}
{\pscal{\xi_{kT}}{\Xbar(y_{kT})^{-1}\xi_{kT}}\geqs h^2, \tau_{\cD_0}\geqs kT}.
\end{align}
Corollary~\ref{cor_timeone} provides a bound for $P_1(h)$, and 
Lemma~\ref{l_endpoint} provides a bound for $Q_1(h)$. Subsequent bounds are
computed by induction, and the following proposition describes one induction
step.

\begin{prop}    \label{prop_induction}
Let $\hat\kappa \leqs \kappa^+(0)\wedge\kappa^\prime$. 
Assume that for some $k\in\N$, 
\begin{align}
\label{timetwo9a}
P_k(h) &\leqs D_k \e^{-\hat\kappa h^2/\sigma^2}, \\
\label{timetwo9b}
Q_k(h) &\leqs \widehat D_k \e^{-\hat\kappa h^2/\sigma^2}.
\end{align}
Then the same bounds hold for $k$ replaced by $k+1$, provided 
\begin{align}   
\label{timetwo10a}
D_{k+1} &\geqs D_k + C^+_{n,m,\gamma}(T,\eps)\widehat D_k 
\frac{\gamma}{\gamma - \hat\kappa} \e^{(\gamma-\hat\kappa)h^2/\sigma^2} \\
\label{timetwo10b}
\widehat D_{k+1} &\geqs \widehat D_k + \widehat C_{n,m,\gamma}(T,\eps).
\end{align}
\end{prop}
\begin{remark}   \label{remark_induction}
Below we will optimize with respect to $\hat\kappa$, but note that in
the case $ \kappa^+(0) = \kappa^\prime= \gamma$, we may either choose
$\hat\kappa < \kappa^+(0)\wedge\kappa^\prime$, or we may
replace~\eqref{timetwo10a} by 
\begin{equation}   \label{timetwo10aa}
D_{k+1} \geqs D_k + C^+_{n,m,\gamma}(T,\eps)\widehat D_k 
\biggbrak{1 + \log \biggpar{\frac{C^+_{n,m,\gamma}(T,\eps)}{\widehat
      D_k}\e^{\gamma h^2/\sigma^2}}}. 
\end{equation}
\end{remark}

\begin{proof}[{\sc Proof of Proposition~\ref{prop_induction}}]
We start by establishing~\eqref{timetwo10b}. The Markov property allows for
the decomposition
\begin{align}
\nonumber
& Q_{k+1}(h)\\ 
\nonumber
& {}\leqs \bigprobin{0,(0,0)}{\tau_{\cB(h)}<kT, \tau_{\cD_0}\geqs kT} \\
\nonumber
& {}+ \Bigexpecin{0,(0,0)}{\indexfct{\tau_{\cB(h)}\geqs kT}
\bigprobin{kT,(\xi_{kT},0)}{\pscal{\xi_{(k+1)T}}{\Xbar(y_{(k+1)T})^{-1}
    \xi_{(k+1)T}}\geqs h^2,  \tau_{\cD_0}\geqs (k+1)T}} \\ 
\label{timetwo11c}
& {}\leqs Q_k(h) + \widehat C_{n,m,\gamma}(T,\eps) \e^{-\hat\kappa
  h^2/\sigma^2}, 
\end{align}
where the initial condition $(\xi_{kT},0)$ indicates that at time $kT$, we
also restart the process of the deterministic slow variables $\ydet_t$ in the
point $y_{kT}\in\cD_0$. In the second line, we used
Lemma~\ref{l_endpoint}. This shows~\eqref{timetwo10b}.

As for~\eqref{timetwo10a}, we again start from a decomposition, similar
to~\eqref{timetwo11c}:
\begin{align}
\nonumber
P_{k+1}(h) & {}= \bigprobin{0,(0,0)}{\tau_{\cB(h)}<kT \wedge \tau_{\cD_0}} \\
\label{timetwo11a}
& \phantom{{}={}}{} +  \Bigexpecin{0,(0,0)}{\indexfct{\tau_{\cB(h)}\geqs kT}
\bigprobin{kT,(\xi_{kT},0)}{\tau_{\cB(h)}<(k+1)T \wedge \tau_{\cD_0}}}.
\end{align}
Corollary~\ref{cor_timeone} allows us to estimate 
\begin{align}
\nonumber
& P_{k+1}(h) \\
\nonumber
& \quad\leqs P_k(h) +
\Bigecondin{0,(0,0)}{\indexfct{\pscal{\xi_{kT}}{\Xbar(y_{kT})^{-1}\xi_{kT}} 
\leqs h^2}
\bigbrak{\varphi\bigpar{\pscal{\xi_{kT}}{\Xbar(y_{kT})^{-1}\xi_{kT}}}
  \wedge1}}{\tau_{\cD_0}\geqs kT} \\
\label{timetwo11b}
&\quad\phantom{{}\leqs P_k(h)+{}}
{}\times\bigprobin{0,(0,0)}{\tau_{\cD_0}\geqs kT},  
\end{align}
with
\begin{equation}   \label{timetwo12}
\varphi(s) = C^+_{n,m,\gamma}(T,\eps) \e^{(\gamma-\hat\kappa)h^2/\sigma^2} 
\e^{-\gamma(h^2-s)/\sigma^2}.
\end{equation}
\eqref{timetwo9b} shows that 
\begin{equation}   \label{timetwo13a}
\bigpcondin{0,(0,0)}
{\pscal{\xi_{kT}}{\Xbar(y_{kT})^{-1}\xi_{kT}}< s}{\tau_{\cD_0}\geqs kT}
\geqs \varphi_k(s),
\end{equation}
where
\begin{equation}   \label{timetwo13b}
\varphi_k(s) \defby \bigpar{1 - \widehat D_k \e^{-\hat\kappa
    s/\sigma^2}}\big/\bigprobin{0,(0,0)}{\tau_{\cD_0}\geqs kT}.
\end{equation}
The functions $\varphi$ and $\varphi_k$ fulfil the assumptions of
Lemma~\ref{l_intbyparts} with 
\begin{equation}   \label{timetwo14}
\e^{\gamma s_2/\sigma^2} = C^+_{n,m,\gamma}(T,\eps)^{-1} \e^{\hat\kappa
h^2/\sigma^2} 
\qquad\text{and}\qquad
\e^{\hat\kappa s_1/\sigma^2} = \widehat D_k.
\end{equation}
For $h^2\leqs s_1$, \eqref{timetwo9a} becomes trivial, while for $h^2>s_1$,
Lemma~\ref{l_intbyparts} shows 
\begin{align}
\nonumber
P_{k+1}(h) 
& \leqs P_k(h) - \varphi(h^2 \wedge s_2)\bigbrak{1-\bigprobin{0,(0,0)}
{\pscal{\xi_{kT}}{\Xbar(y_{kT})^{-1}\xi_{kT}}< h^2, \tau_{\cD_0}\geqs kT}} \\ 
\nonumber
&\phantom{{}\leqs{}} 
{}+ \varphi(s_1) 
+ \int_{s_1}^{s_2 \wedge h^2} \varphi^\prime(s) \widehat D_k
\e^{-\hat\kappa s/\sigma^2} \6s \\
\label{timetwo15}
& \leqs P_k(h) 
+ C^+_{n,m,\gamma}(T,\eps) \widehat D_k \frac{\gamma}{\gamma - \hat\kappa}
\e^{(\gamma-\hat\kappa)h^2/\sigma^2} \e^{-\hat\kappa h^2/\sigma^2}.
\end{align}
Now, \eqref{timetwo10a} is immediate. 
\end{proof}

\goodbreak

Repeated application of the previous result finally leads to the following
estimate.

\begin{cor}   \label{cor_timetwo}
Assume that $y_0\in\cD_0$, $x_0=\bar x(y_0,\eps)$. Then, for every $t>0$, we
have  
\begin{align}
\nonumber
&\bigprobin{0,(x_0,y_0)}{\tau_{\cB(h)}<t\wedge\tau_{\cD_0}} \\
\label{timetwo16}
&\qquad{}\leqs C^+_{n,m,\gamma}(T,\eps) \biggbrak{1 + \widehat
C_{n,m,\gamma}(T,\eps) \biggpar{\frac12+\frac tT}^2
\frac{\gamma}{2(\gamma-\hat\kappa)}}  \e^{-(2\hat\kappa-\gamma) h^2/\sigma^2}.
\end{align}
In addition, the distribution of the endpoint $\x_t$ satisfies 
\begin{equation}
\label{timetwo17}
\bigprobin{0,(x_0,y_0)}
{\pscal{\xi_t}{\Xbar(y_t)^{-1}\xi_t} \geqs h^2, \tau_{\cD_0}\geqs t}  
\leqs \widehat C_{n,m,\gamma}(T,\eps) \biggintpartplus{\frac tT} 
\e^{-\hat\kappa h^2/\sigma^2}.
\end{equation}
\end{cor}
\begin{proof}
We already know the bounds~\eqref{timetwo9a} and~\eqref{timetwo9b} to hold
for $k=1$, with $D_1 = C^+_{n,m,\gamma}(T,\eps)$ and $\widehat D_1 = \widehat
C_{n,m,\gamma}(T,\eps)$. Now the inductive relations~\eqref{timetwo10a}
and~\eqref{timetwo10b} are seen to be satisfied by  
\begin{align} 
\nonumber
\widehat D_k &= k \widehat C_{n,m,\gamma}(T,\eps), \\
\label{timetwo18}
D_k &= C^+_{n,m,\gamma}(T,\eps) \biggbrak{1+\widehat C_{n,m,\gamma}(T,\eps)
\frac{\gamma}{\gamma-\hat\kappa}  \e^{(\gamma-\hat\kappa) h^2/\sigma^2}
\sum_{j=1}^{k-1}j}.
\end{align}
The conclusion follows by taking $k=\intpartplus{t/T}$ and bounding the sum
by $\tfrac12(t/T)(t/T+1)\leqs \tfrac12(t/T + 1/2)^2$. 
\end{proof}

To complete the proof of Part~(a) of Theorem~\ref{thm1}, we first optimize our
choice of $\hat\kappa$, taking into account the constraint $\hat\kappa \leqs
\kappa^+(0) \wedge \kappa^\prime$. By doing so, we find that
\begin{equation}    \label{timetwo19}
\frac{\gamma}{2(\gamma-\hat\kappa)}  \e^{-(2\hat\kappa-\gamma) h^2/\sigma^2}
\leqs \frac{2h^2}{\sigma^2}  \e^{-\kappa h^2/\sigma^2},
\end{equation}
where we have set
\begin{equation}    \label{timetwo20}
\kappa = \gamma \bigbrak{1 - \Order{h} - \Order{\Delta} - \Order{m\eps\rho^2} -
  \bigOrder{\e^{-\text{\it const}/\eps}/(1-2\gamma)}}. 
\end{equation}
Simplifying the prefactor in~\eqref{timetwo16} finally yields the
upper bound
\begin{equation}   \label{timetwo21}
\begin{split}
&\bigprobin{0,(x_0,y_0)}{\tau_{\cB(h)}<t\wedge\tau_{\cD_0}} \\
&\qquad{}\leqs \text{\it const\ }
\frac{(1+t)^2}{\Delta\eps}
\biggbrak{\frac1{(1-2\gamma)^n}
+ \e^{n/4} + \e^{m/4}}
\biggpar{1+\frac{h^2}{\sigma^2}} \e^{-\kappa h^2/\sigma^2}.
\end{split}
\end{equation}

Note that the lower bound in Part~(b) of Theorem~\ref{thm1} is a direct
consequence of the lower bound in Corollary~\ref{cor_timeone}, so that only
Part~(c) remains to be proved.


\subsection{Approaching the adiabatic manifold}   \label{ssec_approach}

The following result gives a rather rough description of the behaviour of
paths starting at a (sufficiently small) distance of order~$1$ from the
adiabatic manifold. It is, however, sufficient to show that with large
probability, these paths will reach the set $\cB(h)$, for some $h>\sigma$,
in a time of order $\eps\abs{\log h}$. 

\begin{prop}   \label{prop_approach}
Let $t$ satisfy the hypotheses of Corollary~\ref{cor_timeone}. Then there
exist constants $h_0$, $\delta_0$, $c_0$ and $K_0$ such that, for $h\leqs
h_0$, $\delta\leqs\delta_0$, $\gamma\in(0,1/2)$ and $\Delta>0$ sufficiently
small,  
\begin{align}   \label{approach1}
\sup_{\xi_0 \,\colon\,\pscal{\xi_0}{\Xbar(y_0)^{-1}\xi_0} \leqs\delta^2}
& \biggprobin{0,(\xi_0,0)}{
\sup_{0 \leqs s \leqs t\wedge\tau_{\cD_0}} 
\frac{\bigpscal{\xi_s}{\Xbar(y_s)^{-1}\xi_s}}
{(h+c_0\delta\e^{-K_0s/\eps})^2} \geqs 1} \\ 
\nonumber
& \leqs
\biggintpartplus{\frac t{\Delta\eps}} 
\biggbrak{\frac1{(1-2\gamma)^{n/2}}
+ \bigpar{\e^{n/4} + 2\e^{m/4}}\e^{-\kappa h^2/\sigma^2} }
\e^{-\kappa h^2/\sigma^2},
\end{align}
where $\kappa = \gamma\brak{1 - \Order{h} - \Order{\Delta} -
\Order{m\eps\rho^2} - \Order{\delta}}$.
\end{prop}

\begin{proof}
We start again by considering an interval $[s,t]$ with $t-s = \Delta\eps$. Let
$\ydet_0=y_0\in\cD_0$. Then 
\begin{align}
\nonumber
P & =\biggprobin{0,(\xi_0,0)}{\sup_{s \leqs u \leqs
    t\wedge\tau_{\cD_0}\wedge\tau_\eta} 
  \frac{\bigpscal{\xi_u}{\Xbar(y_u)^{-1}\xi_u}}{(h+c_0\delta\e^{-K_0u/\eps})^2}
  \geqs 1}  \\
\label{approach2}
&\leqs \biggprobin{0,(\xi_0,0)}{\sup_{s \leqs u \leqs t\wedge\tau} 
\bigpscal{\xi_u}{\Xbar(\ydet_u)^{-1}\xi_u} \geqs H^2},
\end{align}
where $\tau$ is a stopping time defined by
\begin{equation}
\tau = \tau_{\cD_0}\wedge \tau_\eta \wedge
\inf\bigsetsuch{u\geqs0}{\pscal{\xi_u}{\Xbar(\ydet_u)^{-1}\xi_u}\geqs 
(h+c_0\delta\e^{-K_0u/\eps})^2(1-\Order{\Delta})}, 
\end{equation}
and $H^2$ is a shorthand for $H^2= H^2_t =
(h+c_0\delta\e^{-K_0t/\eps})^2(1-\Order{\Delta})$. 

The probability on the right-hand side of~\eqref{approach2} can be bounded, as
in Proposition~\ref{prop_xi}, by the sum $P_0+P_1+P_2$, defined
in~\eqref{timeone13}, provided $H_0+H_1+H_2 = H$. Since
$\norm{U(s)}$ decreases like $\e^{-K_0s/\eps} =
\e^{K_0\Delta}\e^{-K_0t/\eps}$, we have  
\begin{equation}   \label{approach3}
P_0 \leqs \frac1{(1-2\gamma)^{n/2}} 
\exp\Bigbrak{-\gamma\frac{H_0^2 - \text{{\it
const }}\delta^2\e^{2K_0\Delta}\e^{-2K_0t/\eps}}{\sigma^2}}.
\end{equation}
Following the proof of Lemma~\ref{l_nongaussian}, and taking into account
the new definition of $\tau$, we further obtain that 
\begin{equation}   \label{approach4}
P_1 \leqs \e^{n/4}
\exp\Bigset{-\frac{H_1^2}{\sigma^2}\frac1{M_1^2\text{{\it const }}
\brak{(h+h_1)^2 \Psi(t) + c_0^2\delta^2 (t/\eps) \e^{-2K_0t/\eps}}}}.
\end{equation}
As for $P_2$, it can be estimated trivially, provided 
\begin{equation}   \label{approach5}
H_2 \geqs \text{{\it const }} \frac M{K_0} \Bigbrak{(h^2 + hh_1 +
m\eps\rho^2\sigma^2)\Theta(t) 
+ c_0^2\delta^2\e^{K_0\Delta}\e^{-K_0t/\eps}}.
\end{equation}
Choosing $H_1$ in such a way that the exponent in~\eqref{approach4} equals
$H^2/\sigma^2$, we obtain 
\begin{equation}   \label{approach6}
\begin{split}
P \leqs {}&{} 
\biggpar{\frac1{(1-2\gamma)^{n/2}}
+ \e^{n/4}\e^{-H^2/(2\sigma^2)}} \\
&{}\times
\exp\Bigset{-\gamma\frac{H^2}{\sigma^2} \bigbrak{1 - \Order{\Delta} -
\Order{m\eps\rho^2} - \Order{h+h_1+c_0\delta}}}, 
\end{split}
\end{equation}
where we choose $h_1$ proportional to $h + c_0\delta\e^{K_0\Delta}
\e^{-K_0t/\eps}$. The remainder of the proof is similar to the proofs of
Lemma~\ref{l_eta} and Corollary~\ref{cor_timeone}. 
\end{proof}

The preceding lemma shows that after a time $t_1$ of order 
$\eps\abs{\log h}$, the paths are likely to have reached $\cB(h)$. As in 
Lemma~\ref{l_endpoint}, an improved bound for the distribution of the endpoint
$\xi_{t}$ can be obtained. Repeating the arguments leading to
Part~(a) of Theorem~\ref{thm1}, namely using Lemma~\ref{l_intbyparts} on
integration by parts and mimicking the proof of Corollary~\ref{cor_timetwo},
one can show that after any time $t_2\geqs t_1$, the probability of
leaving $\cB(h)$ behaves as if the process had started on the
adiabatic manifold, i.\,e., 
\begin{equation}  \label{approach7}
\biggprobin{0,(\xi_0,0)}{\sup_{t_2 \leqs s \leqs t\wedge\tau_{\cD_0}} 
\bigpscal{\xi_s}{\Xbar(y_s)^{-1}\xi_s} \geqs h^2}
\leqs \cC^+_{n,m,\gamma,\Delta}(t,\eps)  \biggpar{1 + \frac{h^2}{\sigma^2}} 
\e^{-\kappa^+ h^2/\sigma^2},
\end{equation}
uniformly for all $\x_0$ such that $\pscal{\xi_0}{\Xbar(y_0)^{-1}\xi_0} \leqs
\delta^2$. Here $\cC^+_{n,m,\gamma,\Delta}(t,\eps)$ is the same prefactor as
in Theorem~\ref{thm1}, cf.~\eqref{randst18a.6}, and 
\begin{equation}  \label{approach8}
\kappa^+ = \gamma\bigbrak{1 - \Order{h} - \Order{\Delta} - \Order{m\eps\rho^2} 
- \bigOrder{\delta\e^{-\text{\it const}\,(t_2\wedge1)/\eps}/(1-2\gamma)}}.
\end{equation}
This completes our discussion of general initial conditions and, in
particular, the proof of Theorem~\ref{thm1}.


\section{Proofs -- Dynamics of $\z_t$}

In this section, we consider again the SDE
\begin{equation}   \label{zeta1}
\begin{split}
\6x_t & = \frac1\eps f(x_t,y_t,\eps) \6t + \frac{\sigma}{\sqrt\eps}
F(x_t,y_t,\eps) \6W_t, \\
\6y_t & = g(x_t,y_t,\eps) \6t + \sigma^\prime G(x_t,y_t,\eps) \6W_t
\end{split}
\end{equation}
under Assumption~\ref{assump_sde}, that is, when starting near a uniformly
asymptotically stable manifold. We denote by $(\xdet_t,\ydet_t)$, with
$\xdet_t=\bx(\ydet_t,\eps)$, the deterministic solution starting in
$\ydet_0=y_0\in \cD_0$. The system can be rewritten in the
form~\eqref{exit3}, or, in compact notation, as
\begin{equation}   \label{zeta2}
\6\z_t = \bigbrak{\cA(\ydet_t,\eps)\z_t + \cB(\z_t,t,\eps)} \6t 
+ \sigma \bigbrak{\cF_0(\ydet_t,\eps) + \cF_1(\z_t,t,\eps)} \6W_t, 
\end{equation}
where $\transpose\z = (\transpose\xi,\transpose\eta)$, $\cA$ and $\cF_0$
have been defined in~\eqref{randst25}, and the components of $\transpose\cB
= (\eps^{-1}\transpose b,\transpose c)$ and $\transpose{\cF_1} =
(\eps^{-1/2}\transpose{F_1},\rho\transpose{G_1})$ satisfy the
bounds~\eqref{exit4}.

The solution of~\eqref{zeta2} with initial condition $\transpose\z_0 =
(\transpose\xi_0,0)$ can be written in the form 
\begin{align}
\nonumber
\z_t ={}& \cU(t) \z_0 + \sigma\int_0^t \cU(t,s) \cF_0(\ydet_s,\eps)\6W_s \\
&{}+ \int_0^t \cU(t,s) \cB(\z_s,s,\eps)\6s 
+ \sigma\int_0^t \cU(t,s) \cF_1(\z_s,s,\eps)\6W_s. 
\label{zeta3}
\end{align}
The components of the principal solution $\cU(t,s)$ satisfy the bounds 
\begin{align}
\nonumber
\norm{U(t,s)} &\leqs \text{\it const\,} \e^{-K_0(t-s)/\eps}, \\
\label{zeta4}
\norm{S(t,s)} &\leqs \text{\it const\,} \norm{C}_\infty \frac\eps{K_0} 
\bigpar{1 - \e^{-K_0(t-s)/\eps}} \sup_{s\leqs u\leqs t} \norm{V(t,u)}. 
\end{align}
We want to estimate the first-exit time
\begin{equation}   \label{zeta4.5}
\tau_\z = \inf\bigsetsuch{u\geqs0}{\pscal{\z_u}{\cZbar(u)^{-1}\z_u} \geqs h^2},
\end{equation}
with $\cZbar(u)$ defined in~\eqref{randst33}. The inverse of $\cZbar(u)$ is
given by  
\begin{equation}   \label{zeta5}
\cZbar^{-1} = 
\begin{pmatrix}
\lower10pt\vbox{\kern0pt}(\Xbar - \Zbar Y^{-1} \transpose{\Zbar})^{-1} &
-\Xbar^{-1} \Zbar(Y - \transpose{\Zbar}\Xbar^{-1}\Zbar)^{-1} \\
\raise10pt\vbox{\kern0pt} -Y^{-1} \transpose{\Zbar}(\Xbar -  \Zbar Y^{-1}
\transpose{\Zbar})^{-1} & 
(Y - \transpose{\Zbar}\Xbar^{-1}\Zbar)^{-1} 
\end{pmatrix}.
\end{equation}
Since we assume $\norm{\Xbar}_\infty$ and $\norm{\Xbar^{-1}}_\infty$ to be
bounded, $\norm{\Zbar}_\infty = \Order{\sqrt\eps\rho + \eps}$ and 
$\norm{Y^{-1}}_{[0,t]} = \Order{1/(\rho^2+\eps)}$, we have 
\begin{equation}   \label{zeta6}
\cZbar^{-1} = 
\begin{pmatrix}
\lower10pt\vbox{\kern0pt}\Order{1} & \Order{1} \\
\raise10pt\vbox{\kern0pt}\Order{1} & \Order{1/(\rho^2+\eps)}
\end{pmatrix}.
\end{equation}
As in Section~\ref{sec_exit}, we start by examining the dynamics of $\z_u$
on an interval $[s,t]$ with $\Delta = (t-s)/\eps = \orderone{\eps}$. 

The following functions will play a similar r\^ole as the functions $\Phi$
and $\Psi$, introduced in~\eqref{timeone1}, played in
Section~\ref{sec_exit}: 
\begin{align}
\nonumber
\widehat{\Phi}(t) &= \int_0^t \Tr\bigbrak{\transpose{\cJ(v)}
  \transpose{\cU(t,v)} 
\cZbar(t)^{-1} \cU(t,v) \cJ(v)} \6v, \\
\label{zeta6a}
\widehat{\Psi}(t) &= \int_0^t \bignorm{\transpose{\cJ(v)} \transpose{\cU(t,v)}
\cZbar(t)^{-1} \cU(t,v) \cJ(v)} \6v, 
\end{align}
where 
\begin{equation}   \label{zeta6b}
\cJ(v) = \frac1{\sqrt2 \mskip1.5mu M_1 h \norm{\cZbar}_\infty^{1/2}}
\cF_1(\z_v,v,\eps)
= \begin{pmatrix}
\lower10pt\vbox{\kern0pt}\Order{\frac1{\sqrt\eps}} \\ 
\raise10pt\vbox{\kern0pt}\Order{\rho}
\end{pmatrix}
\end{equation}
for $v\leqs\tau_\z$. Using the representations~\eqref{randst27} of $\cU$
and~\eqref{zeta5} of $\cZbar^{-1}$ and expanding the matrix product, one
obtains the relations  
\begin{align}
\nonumber
\widehat{\Phi}(t) &\leqs \Phi(t) + \rho^2 \int_0^t
\Tr\bigbrak{\transpose{V(t,v)} Y(t)^{-1} V(t,v)} \6v 
+ \bigOrder{(n+m)(1+\chi^{(1)}(t)+\chi^{(2)}(t))},\\
\widehat{\Psi}(t) &\leqs \Psi(t) + \rho^2 \int_0^t
\bignorm{\transpose{V(t,v)} Y(t)^{-1} V(t,v)} \6v 
+ \bigOrder{1+\chi^{(1)}(t)+\chi^{(2)}(t)},
\label{zeta6c}
\end{align}
valid for all $t\leqs\tau_\z$. Now we are ready to establish the following
analogue of Proposition~\ref{prop_xi}.

\begin{prop}   \label{prop_zeta}
Fix an initial condition $(x_0,y_0)$ with $y_0\in\cD_0$ and
$x_0=\bx(y_0,\eps)$, and let $t$ be such that $\ydet_u\in\cD_0$ for all
$u\leqs t$. Then, for all $\alpha\in[0,1]$, all $\gamma\in(0,1/2)$ and all
$\mu>0$,  
\begin{align}   \label{zeta7}
\nonumber
&\sup_{\z_0=(\xi_0,0) \,\colon\,\pscal{\xi_0}{\Xbar(y_0)^{-1}\xi_0} 
  \leqs \alpha^2h^2
  }\biggprobin{0,\z_0}{\sup_{s \leqs u \leqs t\wedge\tau_{\cD_0}} 
  \bigpscal{\z_u}{\cZbar(u)^{-1}\z_u} \geqs h^2} \\
\nonumber
& \quad
\leqs \frac{\e^{\Order{m\eps\rho^2}}}{(1-2\gamma)^{(n+m)/2}} \\
\nonumber
& \quad
\hphantom{\leqs{}}
{}\times
\exp\biggset{-\gamma\frac{h^2}{\sigma^2}
  \Bigbrak{1-\alpha^2 - \bigOrder{\Delta + \eps + \mu h 
  + h \norm{\cZbar}_{[0,t]}\bigpar{1+\norm{Y^{-1}}_{[0,t]}^{1/2}\chi^{(1)}(t)}}}} \\
& \quad
\hphantom{\leqs{}}
{}+
\e^{\widehat\Phi(t)/4\widehat\Psi(t)}
\exp{\biggset{-\frac{h^2}{\sigma^2}
\frac{\mu^2(1-\Order{\Delta})}{16M_1^2\norm{\cZbar}_{[0,t]}\widehat\Psi(t)}}}
\end{align}
holds whenever $\sqrt{m\sigma^2}\leqs h < 1/\mu$.
\end{prop}
\begin{proof}
Writing $\z_u = \cU(u,s)\Upsilon_u$, we have 
\begin{equation}   \label{zeta8}
\biggprobin{0,\z_0}{\sup_{s \leqs u \leqs t\wedge\tau_{\cD_0}} 
  \bigpscal{\z_u}{\cZbar(u)^{-1}\z_u} \geqs h^2}
= \biggprobin{0,\z_0}{\sup_{s \leqs u \leqs t\wedge\tau_{\cD_0}\wedge \tau_\z}
  \norm{\cQ(u)\Upsilon_u} \geqs h}
\end{equation}
where $\cQ(u)$ is the symmetric matrix defined by 
\begin{equation}   \label{zeta9}
\cQ(u)^2 = \transpose{\cU(u,s)}\cZbar(u)^{-1}\cU(u,s).
\end{equation}
As in the proof of Proposition~\ref{prop_xi}, we want to eliminate the
$u$-dependence of $\cQ$ in~\eqref{zeta8}. It turns out that the relation
$\norm{\cQ(u)\cQ(t)^{-1}} = 1 + \Order{\Delta}$ still holds in the present
situation, although the proof is less straightforward than before. We
establish this result in Lemma~\ref{lem_Q} below. 

Splitting $\Upsilon_u$ into the sum $\Upsilon_u = \Upsilon^0_u
+ \Upsilon^1_u + \Upsilon^2_u$, where the $\Upsilon^i_u$ are defined in a
way analogous to~\eqref{timeone12}, we can estimate the probability
in~\eqref{zeta8} by the sum $P_0 + P_1 + P_2$, where 
\begin{align}   
\nonumber
P_0 & = \biggprobin{0,\z_0}{\sup_{s \leqs u\leqs t\wedge\tau_{\cD_0}}
  \norm{\cQ(t) \Upsilon^0_u}\geqs H_0}, \\
\label{zeta10}
P_1 & = \biggprobin{0,\z_0}{\sup_{s \leqs u\leqs
    t\wedge\tau_{\cD_0}\wedge\tau_\z} \norm{\cQ(t)\Upsilon^1_u}\geqs H_1}, \\ 
\nonumber
P_2 & =  \biggprobin{0,\z_0}{\sup_{s \leqs u\leqs
    t\wedge\tau_{\cD_0}\wedge\tau_\z} \norm{\cQ(t)\Upsilon^2_u}\geqs H_2}, 
\end{align} 
and $H_0 + H_1 + H_2 = h(1-\Order{\Delta})$. Following the proof of
Lemma~\ref{l_gaussian}, it is straightforward to show that 
\begin{equation}   \label{zeta11}
P_0 \leqs \frac1{(1-2\gamma)^{(n+m)/2}} \exp\Bigset{-\frac\gamma{\sigma^2}
(H_0^2 - \alpha^2 h^2)(1 - \Order{\eps})},
\end{equation}
the sole difference being the factor $\Order{\eps}$ in the exponent
which stems from the fact that $\pscal{\z_0}{\cZbar(0)^{-1}\z_0} =
\pscal{\xi_0}{\Xbar(0)^{-1}\xi_0}(1 + \Order{\eps})$.
Furthermore, similar arguments as in the proof of Lemma~\ref{l_nongaussian}
lead to the bound
\begin{equation}   \label{zeta12}
P_1 \leqs \exp\biggset{- \frac{(H_1^2 - 2\sigma^2M_1^2h^2\norm{\cZbar}_{[0,t]}
\widehat\Phi(t))^2}{16\sigma^2M_1^2h^2H_1^2\norm{\cZbar}_{[0,t]}
\widehat\Psi(t)}}. 
\end{equation}
Finally, the estimate
\begin{align}
\nonumber
\norm{\cQ(t)\Upsilon^2_{u\wedge\tau_\z}}^2
&\leqs \int_0^{u\wedge\tau_\z} \int_0^{u\wedge\tau_\z} 
\bignorm{\transpose{\cB(\z_v,v,\eps)} \transpose{\cU(t,v)} \cZbar(t)^{-1} 
\cU(t,w) \cB(\z_w,w,\eps)} \6v\6w \\
\label{zeta13}
&\leqs \text{{\it const }} \Bigbrak{h^4 \norm{\cZbar}_{[0,t]}^2 
\bigpar{1 + \norm{Y^{-1}}_{[0,t]}\chi^{(1)}(u)^2} +
\bigpar{m\eps\rho^2\sigma^2}^2 \bigpar{1+\chi^{(1)}(u)}}, 
\end{align}
which holds whenever $h\geqs \sqrt{m\sigma^2}$, shows that
$P_2=0$ for 
\begin{equation}   \label{zeta12.5}
H_2 \geqs
\BigOrder{h^2\norm{\cZbar}_{[0,t]}
\bigpar{1+\norm{Y^{-1}}_{[0,t]}^{1/2}\chi^{(1)}(t)} + 
  m\eps\rho^2\sigma^2 \sqrt{1+\chi^{(1)}(t)}}.
\end{equation}
Hence~\eqref{zeta7} follows by taking $H_1 = \mu h^2(1-\Order{\Delta})$.  
\end{proof}

In the proof of Proposition~\ref{prop_zeta}, we have used the following
estimate.

\begin{lemma}   \label{lem_Q}
For $\Delta=(t-s)/\eps$ sufficiently small, 
\begin{equation}   \label{zeta99}
\sup_{s \leqs u \leqs t} \norm{\cQ(u)\cQ(t)^{-1}} = 1 + \Order{\Delta}. 
\end{equation}
\end{lemma}
\begin{proof}
Using the fact that $\cQ(v)^{-2}$ satisfies the ODE
\begin{equation}   \label{zeta99a}
\dtot{}{v} \cQ(v)^{-2} =
\cU(s,v)\cF_0(\ydet_v,\eps)\transpose{\cF_0(\ydet_v,\eps)}
\transpose{\cU(s,v)},
\end{equation}
we obtain the relation
\begin{equation}   \label{zeta99b}
\cQ(u)^2\cQ(t)^{-2} = \one + \cQ(u)^2 \int_u^t
\cU(s,v)\cF_0(\ydet_v,\eps)\transpose{\cF_0(\ydet_v,\eps)}
\transpose{\cU(s,v)}\6v. 
\end{equation}
The definition of $\cF_0$ and the bound~\eqref{zeta4} on $\norm{S}$ allow us
to write
\begin{equation}   \label{zeta99c}
\cU(s,v)\cF_0(\ydet_v,\eps)\transpose{\cF_0(\ydet_v,\eps)}
\transpose{\cU(s,v)} = 
\begin{pmatrix}
\lower10pt\vbox{\kern0pt}\Order{1/\eps} & \Order{\Delta+\rho/\sqrt\eps} \\
\raise10pt\vbox{\kern0pt}\Order{\Delta+\rho/\sqrt\eps} & \Order{\Delta^2\eps +
  \rho^2} 
\end{pmatrix}.
\end{equation}
Using the estimate~\eqref{zeta6} for $\cZbar^{-1}$ and the fact that we
integrate over an interval of length $\Delta\eps$, it follows that 
\begin{equation}   \label{zeta99d}
\cQ(u)^2\cQ(t)^{-2} - \one = \Delta
\begin{pmatrix}
\lower10pt\vbox{\kern0pt}\Order{1} & \Order{\Delta\eps+\rho\sqrt\eps} \\
\raise10pt\vbox{\kern0pt}\Order{1} & \Order{\eps+\rho\sqrt\eps} 
\end{pmatrix},
\end{equation}
which implies~\eqref{zeta99}. 
\end{proof}

Now, Theorem~\ref{thm2} follows from Proposition~\ref{prop_zeta}, by taking a
regular partition of $[0,t]$ with spacing $\Delta\eps$ and
$\mu = 4M_1\smash{\norm{\cZbar}_{[0,t]}^{1/2}}\widehat\Psi(t)^{1/2}$.
We use in particular the fact that $\widehat\Psi(t) = \Order{1 + \chi^{(1)}(t)
  + \chi^{(2)}(t)}$, and that the right-hand side of~\eqref{randst36}
exceeds~$1$ for $h<\sqrt{m\sigma^2}$.


\section{Proofs -- Bifurcations}    \label{sec_bifurcations}

We consider in this section the behaviour of the SDE~\eqref{rp1} near a
bifurcation point. The system can be written in the form 
\begin{equation}   \label{pbif1}
\begin{split}
\6\xi^-_t &= \frac1\eps \hat f^-(\xi^-_t,z_t,y_t,t,\eps) \6t + 
\frac{\sigma}{\sqrt\eps} \widehat F^-(\xi^-_t,z_t,y_t,\eps) \6W_t, \\
\6z_t & = \frac1\eps \hat f^0(\xi^-_t,z_t,y_t,\eps) \6t + 
\frac{\sigma}{\sqrt\eps} \widehat F^0(\xi^-_t,z_t,y_t,\eps) \6W_t, \\
\6y_t & = \hat g(\xi^-_t,z_t,y_t,\eps) \6t + 
\sigma^\prime \widehat G(\xi^-_t,z_t,y_t,\eps) \6W_t,
\end{split}
\end{equation}
compare~\eqref{bif5} and~\eqref{bif6}. We consider the dynamics as long as
$(z_t,y_t)$ evolves in a neighbourhood $\cN$ of the bifurcation point, which is
sufficiently small for the adiabatic manifold to be uniformly asymptotically
stable, that is, all the eigenvalues of $\sdpar {\hat f}x^-(0,z,y,\eps)$
have negative real parts, uniformly bounded away from zero.


\subsection{Exit from $\cB^-(h)$}   \label{ssec_bifexit}

Let $h_\eta, \hz \geqs 0$. In addition to the stopping time 
\begin{equation}   \label{pbif2.0}
\tau_\eta = \inf\bigsetsuch{s>0}{\norm{\eta_s}\geqs h_\eta},
\end{equation}
cf.~\eqref{tau_eta}, we introduce the corresponding stopping time for
$z_s-\zdet_s$, namely, 
\begin{equation}   \label{pbif2}
\tau_z = \inf\bigsetsuch{s>0}{\norm{z_s-\zdet_s}\geqs\hz}.
\end{equation}
The following result is obtained using almost the same line of thought as in
Section~\ref{ssec_timeone}. 

\begin{prop}   \label{prop_pbif1}
Let $t$ be of order $1$ at most. Then, for all initial conditions $\xi^-_0$
such that $\pscal{\xi^-_0}{\Xbarm(y_0,z_0)^{-1}\xi^-_0} \leqs \alpha^2 h^2$
with an $\alpha\in(0,1]$, all $\gamma\in(0,1/2)$, and all sufficiently small
$\Delta>0$,  
\begin{align}  
\nonumber
& \Bigprobin{0,(\xi^-_0,z_0,y_0)}
{\sup_{0\leqs s\leqs t\wedge\tau_\cN\wedge\tau_\eta\wedge\tau_z}
\bigpscal{\xi^-_s}{\Xbarm(y_s,z_s)^{-1}\xi^-_s} \geqs h^2} \\
\nonumber
& \qquad\qquad\quad 
\leqs \biggintpartplus{\frac t{\Delta\eps}} 
\frac{\e^{\Order{(m+q)^{3/2}\sigma}}}{(1-2\gamma)^{(n-q)/2}} 
\exp\biggset{-\gamma\frac{h^2}{\sigma^2}
  \Bigbrak{1-\alpha^2 - \bigOrder{\Delta + (1+\mu) h+h_\eta+h_z}}} \\
\label{pbif3}
& \qquad\qquad\quad \hphantom{\leqs{}}
{} + \biggintpartplus{\frac t{\Delta\eps}} \e^{(n-q)/4}
\exp{\biggset{-\frac{h^2}{\sigma^2}
\frac{\mu^2(1-\Order{\Delta})}{\Order{(1+(h_\eta+h_z)/h)^2}}}}.
\end{align}
\end{prop}
\begin{proof}
The proof is similar to the proof of Corollary~\ref{cor_timeone}, the main
difference being the need for the additional stopping time $\tau_z$. Note that
this results in error terms depending on $h_\eta+h_z$ instead of $h_\eta$
only. For $h^2\geqs (m+q)\sigma^2$, the term $\sigma^2r^- =
\Order{(m+q)\sigma^2}$ yields an error term of order~$h$ in the exponent,
while for $h^2 < (m+q)\sigma^2$, it produces the prefactor
$\e^{\Order{(m+q)^{3/2}\sigma}}$. 
\end{proof}

Next, we need to control the stopping times $\tau_\eta$ and
$\tau_z$. Lemma~\ref{l_eta} holds with minor changes, incorporating the
$z_t$-dependent terms. We find that

\begin{lemma} \label{l_tauetanew}
Let $\xi^-_0$ satisfy $\pscal{\xi^-_0}{\Xbarm(y_0,z_0)^{-1}\xi^-_0}
\leqs h^2$. Then
\begin{align}
\nonumber
&\biggprobin{0,(\xi_0,0)} {\sup_{0 \leqs u \leqs t \wedge
    \tau_{\cB^-(h)}\wedge\tau_z} \norm{\eta_u} \geqs h_\eta} \\
\label{pbif5.5}
&\leqs 2\biggintpartplus{\frac{t}{\Delta\eps}}\e^{m/4} \\
\nonumber
&\hphantom{{}\leqs{}} \times
\exp\biggset{-\kappa_0
\frac{h_\eta^2(1-\Order{\Delta\eps})}{\sigma^2(\rho^2+\eps)\chi^{(2)}(t)}
\biggbrak{1-\biggOrder{\chi^{(1)}(t)\,h_\eta
\Bigpar{1 + \frac{h^2}{h_\eta^2} + \frac{h_z^2}{h_\eta^2} +
(m+q)\frac{\sigma^2}{h_\eta^2}}}}}.  
\end{align}
\end{lemma}
The contribution of $\sigma^2/h_\eta^2$ to the error term might be
puzzling at first glance, but we will apply the preceding lemma for
$h_\eta$ chosen proportional to $h\gg\sigma$, so that
$\sigma^2/h_\eta^2$ will actually be negligible.

The next result allows to control the stopping time $\tau_z$. Let
$U^0(t,s)$ denote the principal solution of $\eps\dot\z =
A^0(\zdet_t,\ydet_t,\eps)\z$, where
$A^0(z,y,\eps)=\partial_zf^0(z,y,\eps)$, and define  
\begin{align}   
\label{pbif6a}
\chi_z^{(1)}(t) 
& = \sup_{0 \leqs s\leqs t} \int_0^s \Bigpar{\sup_{u\leqs v\leqs s}
  \norm{U^0(s,v)}} \6u,  \\ 
\label{pbif6b}
\chi_z^{(2)}(t) 
& = \sup_{0 \leqs s\leqs t} \int_0^s \Bigpar{\sup_{u\leqs v\leqs s} 
\norm{U^0(s,v)}^2} \6u.
\end{align}

\begin{lemma}   \label{l_tauz}
Let $\xi^-_0$ satisfy $\pscal{\xi^-_0}{\Xbarm(y_0,z_0)^{-1}\xi^-_0}
\leqs h^2$. Then
\begin{align}
\nonumber
& \Bigprobin{0,(\xi^-_0,z_0,y_0)}
{\sup_{0\leqs s\leqs t\wedge\tau_{\cB^-(h)}\wedge\tau_\eta}
\norm{z_s - \zdet_s} \geqs \hz} \\
\nonumber
& \qquad
\leqs 2\biggintpartplus{\frac{t}{\Delta\eps}} \e^{q/4} \\
\label{pbif7}
& \qquad\hphantom{{}\leqs{}} \times
\exp\biggset{-\kappa_0
\frac{\eps h_z^2(1-\Order{\Delta\eps})}{\sigma^2 \chi_z^{(2)}(t)}
\biggbrak{1-\biggOrder{\chi_z^{(1)}(t)\,h_z
\Bigpar{1 + \frac{h^2}{h_z^2} + \frac{h_\eta^2}{h_z^2} +
(m+q)\frac{\sigma^2}{h_z^2}}}}}. 
\end{align}
\end{lemma}
\begin{proof}
The proof is almost identical with the proof of Lemma~\ref{l_eta} and
Lemma~\ref{l_tauetanew}, with $\sigma^\prime$ replaced by $\sigma/\sqrt\eps$
and $V$ replaced by $U^0$.  
\end{proof}

Below, we will choose $h_z$ proportional to $h/\sqrt\eps$ for
$h\gg\sigma$, so that the term $(m+q)\sigma^2/h_z^2$ becomes
negligible. 

\begin{proof}[{\sc Proof of Theorem~\ref{thm3}}]
We can repeat the proof of Corollary~\ref{cor_timetwo} in
Section~\ref{ssec_timetwo}, comparing the process to different
deterministic solutions on successive time intervals of length $T$. The
only difference lies in new values for the exponents $\kappa^+(0)$ (resulting
from Proposition~\ref{prop_pbif1}) and $\kappa^\prime$. In fact, choosing
$h_\eta$ proportional to $h$, $h_z$ proportional to
\smash{$(1+\chi_z^{(2)}(T)/\eps)^{1/2} h$} and, finally, $\mu$ proportional to
$1+(h_\eta+h_z)/h$, shows that
\begin{equation}   \label{pbif7.5}
\Bigprobin{0,(\xi^-_0,z_0,y_0)} {\sup_{0\leqs s\leqs T\wedge\tau_\cN}
\bigpscal{\xi^-_s}{\Xbarm(y_s,z_s)^{-1}\xi^-_s} \geqs h^2}
\leqs C_{n,m,q,\gamma}(T,\eps) \e^{-\kappa^+(\alpha) h^2/\sigma^2},
\end{equation}
valid for all $\xi^-_0$ satisfying
$\pscal{\xi^-_0}{\Xbarm(y_0,z_0)^{-1}\xi^-_0} \leqs \alpha^2h^2$ 
and all $T$ of order~$1$ at most. Here
\begin{align}
\label{pbif7.6}
C_{n,m,q,\gamma}(T,\eps) &= \biggintpartplus{\frac T{\Delta\eps}}
\biggbrak{\frac{\e^{\Order{(m+q)^{3/2}\sigma}}}{(1-2\gamma)^{(n-q)/2}}
+ \e^{(n-q)/4} + 2\e^{m/4} + 2\e^{q/4}}, \\
\kappa^+(\alpha) &= \gamma \biggbrak{1 - \alpha^2 - \Order{\Delta} -
\biggOrder{\biggpar{1+\frac{\chi_z^{(2)}(T)}{\eps}} h}}.
\label{pbif7.7}
\end{align}
Similar arguments as in the proof of Lemma~\ref{l_endpoint} yield a bound of
the form   
\begin{equation}
\label{pbif8}
\Bigprobin{0,(\xi^-_0,z_0,y_0)}
{\bigpscal{\xi^-_T}{\Xbarm(y^{\phantom{-}}_T,
z^{\phantom{-}}_T)^{-1}\xi^-_T} \geqs h^2,\tau_\cN\geqs T } 
\leqs \widehat C \e^{-\kappa^\prime h^2/\sigma^2},
\end{equation}
where
\begin{equation}
\label{pbif9}
\kappa^\prime = \gamma\biggbrak{1 - \Order{\Delta} 
- \biggOrder{\biggpar{1+\frac{\chi_z^{(2)}(T)}{\eps}} h}
-\biggOrder{\frac{\e^{-2K_0T/\eps}}{1-2\gamma}}}. 
\end{equation}
In order for the estimates~\eqref{pbif7.5} and~\eqref{pbif8} to be useful, we
need to take $T$ of order $\eps$. However, this leads to an 
error term of order $1$ in the exponent $\kappa^\prime$, which is due to
the fact that $\xi^-_t$ has too little time to relax to the adiabatic
manifold. In order to find the best compromise, we take $T=\theta\eps
\wedge1$ and optimize over $\theta$. Assume we are in the worst case,
when $\norm{U^0}$ grows exponentially like $\e^{K_+t/\eps}$. Then
$\chi^{(2)}_z(T)$ is of the order $\eps\theta \e^{2K_+\th}$. The choice  
\begin{equation}
\label{pbif10}
\e^{-\theta} = \bigbrak{h(1-2\gamma)}^{1/(2(K_0+K_+))}
\end{equation}
yields an almost optimal error term of order $h^\nu (1-2\gamma)^{1-\nu}
\abs{\log (h(1-2\gamma))}$, with $\nu = K_0/(K_0+K_+)$. The smaller
$K_+$, i.\,e., the slower $\chi_z^{(2)}(t)$ grows, the closer $\nu$ is to one.
\end{proof}


\subsection{The reduced system}   \label{ssec_bifreduc}

Given the SDE~\eqref{pbif1}, we call 
\begin{equation}   \label{pbif11}
\begin{split}
\6z^0_t & = \frac1\eps \hat f^0(0,z^0_t,y^0_t,\eps) \6t + 
\frac{\sigma}{\sqrt\eps} \widehat F^0(0,z^0_t,y^0_t,\eps) \6W_t, \\
\6y^0_t & = \hat g(0,z^0_t,y^0_t,\eps) \6t + 
\sigma^\prime \widehat G(0,z^0_t,y^0_t,\eps) \6W_t
\end{split}
\end{equation}
the \defwd{reduced system} of~\eqref{pbif1}. It is obtained by setting
$\xi^-_t= 0$. Let $\z^0_t=(z^0_t,y^0_t)$ and $\z_t=(z_t-z^0_t,y_t-y^0_t)$.
Subtracting~\eqref{pbif11} from~\eqref{pbif1} and making a Taylor expansion
of the drift coefficient, we find that $(\x^-_t,\z_t)$ obeys the SDE 
\begin{equation}   \label{pbif12}
\begin{split}    
\6\xi^-_t & = \frac1\eps \bigbrak{A^-(\z^0_t,\eps)\xi^-_t +
  b(\xi^-_t,\z_t,\z^0_t,\eps)} \6t  
+ \frac{\sigma}{\sqrt\eps} \widetilde F(\xi^-_t,\z_t,\z^0_t,\eps) \6W_t, \\
\6\z_t & = \frac1\eps \bigbrak{C(\z^0_t,\eps)\xi^-_t + B(\z^0_t,\eps)\z_t +
  c(\xi^-_t,\z_t,\z^0_t,\eps)} \6t 
+ \frac{\sigma}{\sqrt\eps} \widetilde\cG(\xi^-_t,\z_t,\z^0_t,\eps) \6W_t, 
\end{split}
\end{equation}
where $\norm{b}$ is of order $\norm{\x^-}^2 + \norm{\z}^2 +
(m+q)\sigma^2$, $\norm{c}$ is of order $\norm{\x^-}^2 +
\norm{\z}^2$ and $\norm{\widetilde\cG}$ is of order $\norm{\x^-} +
\norm{\z}$, while $\norm{\widetilde F}$ is bounded. The matrices $A^-$, $B$
and $C$ are those defined in~\eqref{bif8}, \eqref{bif18a} and~\eqref{bif18b}. 

For a given continuous sample path $\set{\z^0_t(\w)}_{t\geqs0}$
of~\eqref{pbif12}, we denote by $U_\w$ and $\cV_\w$ the principal
solutions of $\eps\dot\x^- = A^-(\z^0_t(\w),\eps)\x^-$ and $\eps\dot\z =
B(\z^0_t(\w),\eps)\z$. If we further define 
\begin{equation}
\label{pbif13}
\cS_\w(t,s) = \frac1\eps \int_s^t \cV_\w(t,u) C(\z^0_u(\w),\eps)
U_\w(u,s)\6u, 
\end{equation}
we can write the solution of~\eqref{pbif12} as 
\begin{align}
\nonumber
\z_t(\w) ={}& 
\frac\sigma{\sqrt\eps} \int_0^t \cV_\w(t,s)
\widetilde\cG(\xi^-_s(\w),\z_s(\w),\z^0_s(\w),\eps)\6W_s(\w) \\
\nonumber
&{}+ \frac\sigma{\sqrt\eps} \int_0^t \cS_\w(t,s)
\widetilde F(\xi^-_s(\w),\z_s(\w),\z^0_s(\w),\eps)\6W_s(\w) \\
\nonumber
&{}+ \frac1\eps \int_0^t \cV_\w(t,s)
c(\xi^-_s(\w),\z_s(\w),\z^0_s(\w),\eps)\6s \\
\label{pbif14}
&{}+ \frac1\eps \int_0^t \cS_\w(t,s)
b(\xi^-_s(\w),\z_s(\w),\z^0_s(\w),\eps)\6s.
\end{align}
Concerning the first two summands in~\eqref{pbif14}, note that the
identities 
\begin{equation}   \label{pbif14.5}
\begin{split}
\cV_\w(t,s) & = \cV_\w(t,0) \cV_\w(s,0)^{-1}, \\
\cS_\w(t,s) & = \cS_\w(t,0) U_\w(s,0)^{-1} + \cV_\w(t,0) \cS_\w(s,0)^{-1}
\end{split}
\end{equation}
allow to rewrite the stochastic integrals in such a way that the integrands
are adapted with respect to the filtration generated by $\{W_s\}_{s\geqs0}$.

We now assume the existence of a stopping time $\tau\leqs\tau_{\cB^-(h)}$ and
deterministic functions $\vartheta(t,s)$, $\vartheta_C(t,s)$ such that 
\begin{equation}   \label{pbif15}
\begin{split}
\bignorm{\cV_\w(t,s)} &\leqs \vartheta(t,s), \\
\bignorm{\cV_\w(t,s)C(\z^0_s(\w),\eps)} &\leqs \vartheta_C(t,s),
\end{split}
\end{equation}
uniformly in $\eps$, whenever $s\leqs t\leqs\tau(\w)$, and define 
\begin{align}
\label{pbif16a}
\chi^{(i)}(t) &= \sup_{0\leqs s\leqs t}\frac1\eps \int_0^s 
 \vartheta(s,u)^i \6u,  
& i&=1,2,\\
\label{pbif16b}
\chi_C^{(i)}(t) &= \sup_{0\leqs s\leqs t}\frac1\eps \int_0^s 
\Bigpar{\sup_{u\leqs v\leqs s} \vartheta_C(s,v)^i} \6u,
& i&=1,2. 
\end{align}

The following proposition establishes a local version of Theorem~\ref{thm4}.

\begin{prop}
\label{prop_bifreduc}
Let $\Delta$ be sufficiently small, fix times $s<t$ such that
$t-s=\Delta\eps$, and assume that there exists a constant $\vartheta_0>0$ such
that $\vartheta(u,s)\leqs \vartheta_0$ and $\vartheta_C(u,s)\leqs
\vartheta_0$, whenever $u\in[s,t]$. Then there exist constants $\kappa_0,
h_0>0$ such that for all $h\leqs
h_0\brak{\chi^{(1)}(t)\vee\chi_C^{(1)}(t)}^{-1}$, 
\begin{equation}
\label{pbif17}
\Bigprobin{0,0}{\sup_{s\wedge\tau \leqs u < t\wedge\tau} \norm{\z_u}\geqs h} 
\leqs 2\e^{(m+q)/4} \exp\Bigset{-\kappa_0
\frac{h^2}{\sigma^2}\frac1{\chi_C^{(2)}(t) + h \chi_C^{(1)}(t) +
h^2\chi^{(2)}(t)}}.  
\end{equation}
\end{prop}
\begin{proof}
The proof follows along the lines of the proof of Lemma~\ref{l_eta}, the
main difference lying in the fact that the stochastic integrals in
~\eqref{pbif14} involve the principal solutions $U_\w$, $\cV_\w$
depending on the realization of the process. However, the existence of the
deterministic bound~\eqref{pbif15} allows for a similar conclusion. In
particular, the first and second term in~\eqref{pbif14} create respective
contributions of the form
\begin{align}
\label{pbif18a}
&\e^{(m+q)/4} \exp\Bigset{-\frac{H_0^2}
{16\sigma^2 h^2M_1^2\chi^{(2)}(t)}}\\
\label{pbif18b}
&\e^{(m+q)/4} \exp\Bigset{-\frac{H_1^2}
{16\sigma^2 M_1^2\chi_C^{(2)}(t)}}
\end{align}
to the probability~\eqref{pbif17}. The third and fourth term only cause
corrections of order $h\chi^{(1)}(t)$ and
$h\chi_C^{(1)}(t)[1+(m+q)\sigma^2/h^2]$ in the exponent. Note that we
may assume $h^2\gg(m+q)\sigma^2$ as well as
$h\gg(m+q)\sigma^2\chi_C^{(1)}(t)$, because Estimate~\eqref{pbif17} is
trivial otherwise.
\end{proof}

Now Theorem~\ref{thm4} follows from Proposition~\ref{prop_bifreduc} by
using a partition of the interval $[0,t]$ into smaller intervals of length
$\Delta\eps$.


\small
\bibliography{geom}

\begin{thebibliography}{10}

\bibitem{Arnold1}
L.~Arnold.
\newblock {\em Random Dynamical Systems}.
\newblock Springer-Verlag, Berlin, 1998.

\bibitem{Azencott}
R.~Azencott.
\newblock Petites perturbations al\'eatoires des syst\`emes dynamiques:
  d\'eveloppements asymptotiques.
\newblock {\em Bull.\ Sci.\ Math.\ (2)}, 109:253--308, 1985.

\bibitem{Bellman}
R.~Bellman.
\newblock {\em Introduction to Matrix Analysis}.
\newblock McGraw--Hill, New York, 1960.

\bibitem{Benoit}
E.~Beno\^{\i}t, editor.
\newblock {\em Dynamic Bifurcations}, Berlin, 1991. Springer-Verlag.

\bibitem{BG2}
N.~Berglund and B.~Gentz.
\newblock A sample-paths approach to noise-induced synchronization:
  {S}tochastic resonance in a double-well potential.
\newblock To appear in {Ann.}\ {Appl.}\ {Probab.}\\ Available at {\tt
  http://arXiv.org/abs/math.PR/0012267}, 2000.

\bibitem{BG3}
N.~Berglund and B.~Gentz.
\newblock The effect of additive noise on dynamical hysteresis.
\newblock {\em Nonlinearity}, 15(3):605--632, 2002.
\newblock DOI 10.1088/0951-7715/15/3/305.

\bibitem{BG1}
N.~Berglund and B.~Gentz.
\newblock Pathwise description of dynamic pitchfork bifurcations with additive
  noise.
\newblock {\em {Probab.} {Theory} {Related} {Fields}}, 122(3):341--388, 2002.
\newblock DOI 10.1007/s004400100174.

\bibitem{BEGK}
A.~Bovier, M.~Eckhoff, V.~Gayrard, and M.~Klein.
\newblock Metastability in reversible diffusion processes {I}. {S}harp
  asymptotics for capacities and exit times.
\newblock Preprint WIAS-767, 2002.

\bibitem{BGK}
A.~Bovier, V.~Gayrard, and M.~Klein.
\newblock Metastability in reversible diffusion processes {II}. {P}recise
  asymptotics for small eigenvalues.
\newblock Preprint WIAS-768, 2002.

\bibitem{CF1}
H.~Crauel and F.~Flandoli.
\newblock Attractors for random dynamical systems.
\newblock {\em Probab. Theory Related Fields}, 100(3):365--393, 1994.

\bibitem{Day1}
M.~V. Day.
\newblock On the exponential exit law in the small parameter exit problem.
\newblock {\em Stochastics}, 8:297--323, 1983.

\bibitem{Day2}
M.~V. Day.
\newblock On the exit law from saddle points.
\newblock {\em Stochastic Process.\ Appl.}, 60:287--311, 1995.

\bibitem{Fenichel}
N.~Fenichel.
\newblock Geometric singular perturbation theory for ordinary differential
  equations.
\newblock {\em J. Differential Equations}, 31(1):53--98, 1979.

\bibitem{FJ}
W.~H. Fleming and M.~R. James.
\newblock Asymptotic series and exit time probabilities.
\newblock {\em Ann. Probab.}, 20(3):1369--1384, 1992.

\bibitem{Freidlin2}
M.~I. Freidlin.
\newblock On stable oscillations and equilibriums induced by small noise.
\newblock {\em J.~Stat.\ Phys.}, 103:283--300, 2001.

\bibitem{FW}
M.~I. Freidlin and A.~D. Wentzell.
\newblock {\em Random Perturbations of Dynamical Systems}.
\newblock Springer-Verlag, New York, second edition, 1998.

\bibitem{Grad}
I.~S. Grad{\v{s}}te{\u\i}n.
\newblock Application of {A}. {M}. {L}yapunov's theory of stability to the
  theory of differential equations with small coefficients in the derivatives.
\newblock {\em Mat. Sbornik N. S.}, 32(74):263--286, 1953.

\bibitem{Haberman}
R.~Haberman.
\newblock Slowly varying jump and transition phenomena associated with
  algebraic bifurcation problems.
\newblock {\em SIAM J. Appl. Math.}, 37(1):69--106, 1979.

\bibitem{Jones}
C.~K. R.~T. Jones.
\newblock Geometric singular perturbation theory.
\newblock In {\em Dynamical systems (Montecatini Terme, 1994)}, pages 44--118.
  Springer, Berlin, 1995.

\bibitem{Kifer}
Y.~Kifer.
\newblock The exit problem for small random perturbations of dynamical systems
  with a hyperbolic fixed point.
\newblock {\em Israel J. Math.}, 40(1):74--96, 1981.

\bibitem{LebovitzSchaar2}
N.~R. Lebovitz and R.~J. Schaar.
\newblock Exchange of stabilities in autonomous systems.
\newblock {\em Studies in Appl. Math.}, 54(3):229--260, 1975.

\bibitem{LebovitzSchaar1}
N.~R. Lebovitz and R.~J. Schaar.
\newblock Exchange of stabilities in autonomous systems. {I}{I}. {V}ertical
  bifurcation.
\newblock {\em Studies in Appl. Math.}, 56(1):1--50, 1976/77.

\bibitem{MishchenkoRozov}
E.~F. Mishchenko and N.~K. Rozov.
\newblock {\em Differential equations with small parameters and relaxation
  oscillations}.
\newblock Plenum Press, New York, 1980.

\bibitem{Neishtadt1}
A.~I. Ne{\u\i}shtadt.
\newblock Persistence of stability loss for dynamical bifurcations {I}.
\newblock {\em Differential Equations}, 23:1385--1391, 1987.

\bibitem{Neishtadt2}
A.~I. Ne{\u\i}shtadt.
\newblock Persistence of stability loss for dynamical bifurcations {II}.
\newblock {\em Differential Equations}, 24:171--176, 1988.

\bibitem{Pontryagin}
L.~S. Pontryagin.
\newblock Asymptotic behavior of solutions of systems of differential equations
  with a small parameter in the derivatives of highest order.
\newblock {\em Izv. Akad. Nauk SSSR. Ser. Mat.}, 21:605--626, 1957.

\bibitem{Schmal}
B.~Schmalfu{\ss}.
\newblock Invariant attracting sets of nonlinear stochastic differential
  equations.
\newblock In H.~Langer and V.~Nollau, editors, {\em Markov processes and
  control theory}, volume~54 of {\em Math.\ Res.}, pages 217--228, Berlin,
  1989. Akademie-Verlag.
\newblock {G}au\ss ig, 1988.

\bibitem{Tihonov}
A.~N. Tihonov.
\newblock Systems of differential equations containing small parameters in the
  derivatives.
\newblock {\em Mat.\ Sbornik N. S.}, 31:575--586, 1952.

\end{thebibliography}
\bibliographystyle{abbrv}      

\goodbreak

\bigskip\bigskip\noindent
{\small 
Nils Berglund \\ 
{\sc Department of Mathematics, ETH Z\"urich} \\ 
ETH Zentrum, 8092~Z\"urich, Switzerland \\
{\it and} \\
{\sc FRUMAM, CPT--CNRS Luminy} \\
Case 907, 13288~Marseille Cedex 9, France \\
{\it and} \\
{\sc PHYMAT, Universit\'e de Toulon} \\
B.\,P.~132, 83957~La~Garde Cedex, France \\
{\it E-mail address: }{\tt berglund@cpt.univ-mrs.fr}

\bigskip\noindent
Barbara Gentz \\ 
{\sc Weierstra\ss\ Institute for Applied Analysis and Stochastics} \\
Mohrenstra{\ss}e~39, 10117~Berlin, Germany \\
{\it E-mail address: }{\tt gentz@wias-berlin.de}
}


\end{document}